\documentclass[small]{article}
\usepackage{amsmath,amstext,amsthm,amsfonts}
\usepackage{makeidx}
\usepackage{amssymb,latexsym}
\usepackage{txfonts, pxfonts}
\usepackage[english]{babel}
\usepackage[pdftex]{graphicx}
\numberwithin{equation}{section}
\newtheorem{theorem}{\textbf{Theorem}}[section]
\newtheorem{proposition}{\textbf{Proposition}}[section]
\newtheorem{lemma}{\textbf{Lemma}}[section]
\newtheorem{corollary}{Corollary}[section]
\theoremstyle{definition}
\newtheorem{definition}{\textbf{Definition}}[section]
\theoremstyle{remark}
\newtheorem{remark}{\it{Remark}}[section]
\newenvironment{notation}[1][Notation]{\begin{trivlist}
\item[\hskip \labelsep {\bfseries #1}]}{\end{trivlist}}

\def\e{\epsilon}
\def\ei{\epsilon_i}
\def\R{\mathbb{R}}
\def\Rn{{\mathbb{R}}^n_+}

\def\d{\partial}

\def\fermi{\psi_i:B^+_{\delta}(0)\to M}
\def\fermilinha{\psi_i:B^+_{\delta'}(0)\to M}
\def\a{\alpha}
\def\b{\beta}
\def\l{\lambda}

\def\Bei{B^+_{\delta\ei^{-1}}}
\def\Beilinha{B^+_{\delta'\ei^{-1}}}

\def\Jzi{\left(y^b\partial_b\phi_i+\frac{n-2}{2}\phi_i\right)}
\def\JU{\left(y^b\partial_bU+\frac{n-2}{2}U\right)}

\def\ba{\begin{align}}
\def\ea{\end{align}}
\def\bp{\begin{proof}}
\def\ep{\end{proof}}

\title{A compactness theorem for scalar-flat metrics
on manifolds with boundary}

\author{S\'ergio de Moura Almaraz}

\begin{document}

\title{A compactness theorem for scalar-flat metrics
on manifolds with boundary}

\author{S\'ergio de Moura Almaraz}

\maketitle

\begin{abstract}
Let $(M^n,g)$ be a compact Riemannian manifold with boundary $\d M$. This paper is concerned with the set of scalar-flat metrics which are in the conformal class of $g$ and have $\d M$ as a constant mean curvature hypersurface. We prove that this set is compact for dimensions $n\geq 7$ under the  
generic condition that the trace-free 2nd fundamental form of $\d M$ is nonzero everywhere.   
\end{abstract}

\section{Introduction}

In 1960,  H. Yamabe (\cite{yamabe}) raised the following question: 

\vspace{0.2cm}
\noindent
YAMABE PROBLEM: {\it{Given $(M^n,g)$, a compact Riemannian manifold (without boundary) of dimension $n\geq 3$, is there a Riemannian metric, conformal to $g$, with constant scalar curvature?}}
\vspace{0.1cm}


This question was affirmatively answered after the works of Yamabe himself, N. Trudinger (\cite{trudinger}), T. Aubin (\cite{aubin1}) and R. Schoen (\cite{schoen1}). (See \cite{lee-parker} and \cite{schoen-yau3} for nice surveys on the issue.)

In 1992, J. Escobar (\cite{escobar3}) studied the following Yamabe-type problem, for manifolds with boundary:

\vspace{0.2cm}
\noindent
YAMABE PROBLEM (boundary version): {\it{Given $(M^n,g)$, a compact Riemannian manifold of dimension $n\geq 3$ with boundary, is there a Riemannian metric, conformal to $g$, with zero scalar curvature and  constant boundary mean curvature?}}
\vspace{0.1cm}


In analytical terms, the problem proposed by Escobar corresponds to finding a positive solution to 
\begin{align}\label{eq:u'}
\begin{cases}
L_{g}u=0,&\text{in}\:M,
\\
B_{g}u+Ku^{\frac{n}{n-2}}=0,&\text{on}\:\partial M,
\end{cases}
\end{align}
for some constant $K$, where $L_g=\Delta_g-\frac{n-2}{4(n-1)}R_g$ is the conformal Laplacian and $B_g=\frac{\d}{\d\eta}-\frac{n-2}{2}h_g$. Here, $\Delta_g$ is the Laplace-Beltrami operator, $R_g$ is the scalar curvature, $h_g$ is the mean curvature of $\d M$ and  $\eta$ is the inward unit normal vector to $\d M$. 

The solutions of the equations (\ref{eq:u'}) are the critical points of the functional 
$$
Q(u)=\frac{\int_M|\nabla_gu|^2+\frac{n-2}{4(n-1)}R_gu^2dv_g+\frac{n-2}{2}\int_{\d M}h_gu^2d\sigma_g}
{\left(\int_{\d M}|u|^\frac{2(n-1)}{n-2}d\sigma_g\right)^{\frac{n-2}{n-1}}}\,,
$$
where $dv_g$ and $d\sigma_g$ denote the volume forms of $M$ and $\d M$, respectively. 
In order to prove the existence of solutions to the equations (\ref{eq:u'}), Escobar introduced the conformally invariant Sobolev quotient
$$
Q(M,\d M)=\inf\{Q(u);\:u\in C^1(\bar{M}), u\nequiv 0 \:\text{on}\: \d M\}\,.
\vspace{-0.1cm}
$$
The question of existence of solutions to the equations (\ref{eq:u'}) was studied in \cite{almaraz2}, \cite{chen}, \cite{escobar3}, \cite{escobar4}, \cite{ahmedou-felli1}, \cite{coda1} and \cite{coda2}. 
The regularity of these solutions was established in \cite{cherrier}. 
Conformal metrics of constant scalar curvature and zero boundary mean curvature were studied in \cite{brendle-chen}, \cite{escobar2} (see also \cite{ambrosetti-li-malchiodi} and \cite{han-li}).

In the case of manifolds without boundary, the question of compactness of the full set of solutions to the Yamabe equation was first raised by R. Schoen in a topics course at Stanford University in 1988. A necessary condition is that the manifold $M^n$ is not conformally equivalent to the sphere $S^n$. 
This problem was studied in  \cite{druet1}, \cite{druet2}, \cite{li-zhang}, \cite{li-zhang2}, \cite{li-zhu2}, \cite{marques}, \cite{schoen4} and \cite{schoen-zhang} and was completely solved in a series of three papers: \cite{brendle2}, \cite{brendle-marques} and \cite{khuri-marques-schoen}. In \cite{brendle2}, Brendle discovered the first smooth counterexamples for dimensions $n\geq 52$ (see \cite{berti-malchiodi} for nonsmooth examples). In \cite{khuri-marques-schoen}, Khuri, Marques and Schoen proved compactness for dimensions $3\leq n\leq 24$. Their proof contains both a local and a global aspect. The local aspect involves the vanishing of the Weyl tensor at any blow-up point and the global aspect involves the Positive Mass Theorem. Finally, in \cite{brendle-marques}, Brendle and Marques extended the counterexamples of \cite{brendle2} to the remaining dimensions  $25\leq n\leq 51$. In \cite{li-zhang}, \cite{li-zhang2} and \cite{marques} the authors proved compactness for $n\geq 6$ under the condition that the Weyl tensor is nonzero everywhere.

In the present work we are interested in the compactness of the set of positive solutions to 
\begin{align}\label{main:eq}
\begin{cases}
L_{g}u=0,&\text{in}\:M,
\\
B_{g}u+Ku^{p}=0,&\text{on}\:\partial M,
\end{cases}
\end{align}
where $1< p\leq \frac{n}{n-2}$. 
A necessary condition is that $M$ is not conformally equivalent to $B^n$. 
As stated by Escobar in \cite{escobar3}, $Q(M,\d M)$ is positive, zero or negative if the first eigenvalue $\l_1(B_g)$ of the problem 
\begin{align}\notag
\begin{cases}
L_{g}u=0,&\text{in}\:M,
\\
B_{g}u+\l u=0,&\text{on}\:\partial M
\end{cases}
\end{align}
is positive, zero or negative, respectively. If $\l_1(B_g)<0$, the solution to the equations (\ref{main:eq}) is unique. If $\l_1(B_g)=0$, the equations (\ref{main:eq}) become linear and the solutions are unique up to a multiplication by a positive constant. Hence, the only interesting case is the one when $\l_1(B_g)>0$.

We expect that, as in the case of manifolds without boundary, there should be counterexamples to compactness of the set of solutions to the equations (\ref{main:eq}) in high dimensions. In this work we address the question of whether compactness of these solutions holds generically in any dimension. 

Our first result is the following:
\begin{theorem}\label{compactness:thm'}
Let $(M^n,g)$ be a Riemannian manifold with dimension $n\geq 7$ and boundary $\d M$.  
Assume that $Q(M,\d M)>0$. Let $\{u_i\}$ be a sequence of solutions to the equations (\ref{main:eq}) with $p=p_i\in[1+\gamma_0,\frac{n}{n-2}]$ for any small fixed $\gamma_0>0$. Suppose there is a sequence $\{x_i\}\subset\d M$,  $x_i \to x_0$,  of local maxima points of $u_i|_{\d M}$ such that $u_i(x_i)\to \infty$. Then the trace-free 2nd fundamental form of $\d M$ vanishes at $x_0$.
\end{theorem}

By linear elliptic theory, uniform estimates for the solutions of equations (\ref{main:eq})  imply $C^{k,\a}$-estimates, for some $0<\a<1$. By the Harnack-type inequality of Lemma \ref{Harnack:han-li} (proved in \cite{han-li}), uniform estimates on the boundary $\d M$ imply uniform estimates on $M$.  Hence, an immediate consequence of Theorem \ref{compactness:thm'} is a compactness theorem for Riemannian manifolds of dimension $n\geq 7$ that satisfy the condition that the boundary  trace-free 2nd fundamental form is nonzero everywhere.  More precisely:
\begin{theorem}\label{compactness:thm}
Let $(M^n,g)$ be a Riemannian manifold with dimension $n\geq 7$ and boundary $\d M$.  Suppose $Q(M,\d M)>0$ and that the trace-free 2nd fundamental form of $\d M$ is nonzero everywhere. 
Then, given a small $\gamma_0>0$, there exists $C>0$ such that for any  $p\in\left[1+\gamma_0,\frac{n}{n-2}\right]$ and $u>0$ solution to the equations (\ref{main:eq})
we have
$$C^{-1}\leq u\leq C\:\:\:\:\: \text{and}\:\:\:\:\:\|u\|_{C^{2,\a}(M)}\leq C\,,$$
for some $0<\a<1$.
\end{theorem}

It was pointed out to me by F. Marques that a transversality argument implies that the second fundamental form condition above is generic for $n\geq 4$. In other words, the set of the Riemannian metrics on $M^n$ such that the trace-free second fundamental form of $\d M$ is nonzero everywhere is open and dense in the space of all Riemannian metrics on $M$ for $n\geq 4$.

We should mention that Theorem \ref{compactness:thm} does not use the Positive Mass Theorem, since the proof of Theorem \ref{compactness:thm'} contains only a local argument, based in a Pohozaev-type identity. 

The problem of compactness of solutions to the equations (\ref{main:eq}) was also studied by  V. Felli and M. Ould Ahmedou in the conformally flat case with umbilic boundary (\cite{ahmedou-felli1}) and in the three-dimensional case with umbilic boundary (\cite{ahmedou-felli2}). Other compactness results for similar equations  were obtained by Z. Han and Y. Li in \cite{han-li} and by Z. Djadli, A. Malchiodi and M. Ould Ahmedou in \cite{djadli-malchiodi-ahmedou1} and \cite{djadli-malchiodi-ahmedou2}.

A consequence of Theorem \ref{compactness:thm} is the computation of the total Leray-Schauder degree of all solutions to the equations (\ref{eq:u'}), as in \cite{ahmedou-felli1}, \cite{ahmedou-felli2} and \cite{han-li} (see also \cite{khuri-marques-schoen}). When $\lambda_1(B_g)>0$, we can define a map
$F_p:\bar{\Omega}_{\Lambda}\to C^{2,\a}(M)$ by $F_p(u)=u+T(E(u)u^p)$. Here, $E(u)=\int_M|\nabla_gu|^2+\frac{n-2}{4(n-1)}R_gu^2dv_g+\frac{n-2}{2}\int_{\d M}h_gu^2d\sigma_g$ is the energy of $u$, $T$ is the operator defined  by $T(v)=u$, where $u$ is the unique solution to 
\begin{align}
\begin{cases}\notag
L_{g}u=0,&\text{in}\:M\,,
\\
B_{g}u=v,&\text{on}\:\partial M\,,
\end{cases}
\end{align}
and $\Omega_{\Lambda}=\{u\in C^{2,\a}(M);\,|u|_{C^{2,\a}(M)}<\Lambda, \,u>\Lambda^{-1}\}$. From elliptic theory we know that the map $u\mapsto T(E(u)u^p)$ is compact from $\bar{\Omega}_{\Lambda}$ into $C^{2,\a}(M)$. Hence, $F_p$ is of the form $I+\text{compact}$. If $0\neq F_p(\d\Omega_{\Lambda})$, we may define the Leray-Schauder degree (see \cite{nirenberg}) of $F_p$ in the region $\Omega_{\Lambda}$ with respect to $0\in C^{2,\a}(M)$, denoted by $\text{deg}(F_p,\Omega_{\Lambda}, 0)$. Observe that $F_p(u)=0$ if and only if $u$ is a solution to 
\begin{align}\notag
\begin{cases}
L_{g}u=0,&\text{in}\:M,
\\
B_{g}u+E(u)u^{p}=0,&\text{on}\:\partial M.
\end{cases}
\end{align}
Observe that these equations imply that $\int_{\d M}u^{p+1}d\sigma_g=1$.
By the homotopy invariance of the degree, $\text{deg}(F_p,\Omega_{\Lambda}, 0)$ is constant for all $p\in\left[1,\frac{n}{n-2}\right]$ provided that $0\neq F_p(\d\Omega_{\Lambda})$ for all $p\in\left[1,\frac{n}{n-2}\right]$. In the linear case, when $p=1$, we have $\text{deg}(F_1,\Omega_{\Lambda}, 0)=-1$. This is the content of Lemma 4.2 of \cite{ahmedou-felli1}, which is a modification of the arguments in \cite{han-li}, pp.528-529. Thus, for $\Lambda$ sufficiently large, Theorem \ref{compactness:thm} allow us to calculate the degree for all $p\in\left[1,\frac{n}{n-2}\right]$. Hence, we have:
\begin{theorem}
Let $(M^n,g)$ satisfy the assumptions of Theorem \ref{compactness:thm}. Then, for $\Lambda$ sufficiently large and all $p\in\left[1,\frac{n}{n-2}\right]$, we have $\text{deg}(F_p,\Omega_{\Lambda}, 0)=-1$.
 \end{theorem}

We will now outline the proof of Theorem \ref{compactness:thm'}.  
The strategy of the proof is similar to the one proposed by Schoen in the case of manifolds without boundary.  It is based on finding local obstructions to blow-up by means of a Pohozaev-type identity. 
We suppose that there is a simple blow-up point  for a sequence $\{u_i\}$. We then   approximate the sequence $\{u_i\}$ by the standard Euclidean solution plus a correction term $\phi_i$. The function $\phi_i$ is defined as a solution to a non-homogeneous linear equation. We then use the Pohozaev identity to prove that the boundary trace-free 2nd fundamental form  vanishes at the blow-up point. 
Finally we apply the Pohozaev identity to establish, after rescaling arguments, a sign condition that allows the reduction to the simple blow-up case.

An important part  in our proof is the use of the correction term $\phi_i$ to obtain refined pointwise blow-up estimates. The idea of using a correction term first appeared in \cite{hebey-vaugon} and was significantly 
improved in \cite{brendle3}. This type of blow-up estimate was derived in \cite{khuri-marques-schoen} where the authors studied compactness in the case of manifolds without boundary. Although we do not have the kind of explicit control of the terms $\phi_i$ the authors had in \cite{khuri-marques-schoen}, a key observation is that some orthogonality conditions are sufficient to obtain the vanishing of the boundary trace-free 2nd fundamental form.

In Section 2 we establish some notation and discuss some basic results. In Section 3 we prove the Pohozaev identity we will work with.
In Section 4 we discuss the concepts of isolated and isolated simple blow-up points and state some basic properties.
In Section 5 we find the correction term $\phi_i$ and prove its properties.  
In Section 6 we obtain the pointwise estimates for $u_i$. 
In Section 7 we prove the vanishing of the trace-free 2nd fundamental form at any isolated simple blow-up point and prove the Pohozaev sign condition. 
In Section 8 we reduce our analysis to the case of isolated simple blow-up points and prove Theorem \ref{compactness:thm'}. 
\\\\
{\bf{Acknowledgements.}}
The content of this paper is a part of the author's doctoral thesis 
(\cite{almaraz1}). 
The author would like to express his gratitude to his advisor Prof. Fernando C. Marques for numerous mathematical conversations and constant 
encouragement. While the author was at IMPA, he was fully supported by CNPq-Brazil.

\section{Preliminaries}

\subsection{Notations}\label{subsec:notations}
Throughout this work we will make use of the index notation for tensors, commas denoting covariant differentiation. We will adopt the summation convention whenever confusion is not possible. When dealing with coordinates on manifolds with boundary, we will use indices $1\leq i,j,k,l\leq n-1$ and $1\leq a,b,c,d\leq n$. In this context, lines under or over an object mean the restriction of the metric to the boundary is involved.

We will denote by $g$ the Riemannian metric and set $\det g=\det g_{ab}$. The induced metric on $\d M$ will be denoted by $\bar{g}$. We will denote by $\nabla_g$  the covariant derivative and by $\Delta_g$ the Laplacian-Beltrami operator. The full curvature tensor will be denoted by $R_{abcd}$,  the Ricci tensor by $R_{ab}$ and the scalar curvature by $R_g$ or $R$. The second fundamental form of the boundary will be denoted by $h_{kl}$ and the mean curvature,  $\frac{1}{n-1}tr (h_{kl})$, by $h_g$ or $h$.  By $\pi_{kl}$ we will denote the trace-free second fundamental form, $h_{kl}-h\bar{g}_{kl}$. 

By $\Rn$ we will denote the half-space $\{z=(z_1,...,z_n)\in \R^n;\:z_n\geq 0\}$. If $z\in\Rn$ we set $\bar{z}=(z_1,...,z_{n-1})\in\R^{n-1}\cong \d\Rn$. 
We define $B^+_{\delta}(0)=\{z\in\Rn\,;\:|z|<\delta \}$. We also denote $B^+_{\delta}=B^+_{\delta}(0)$ for short. 
We set $\d^+B^+_{\delta}(0)=\d B^+_{\delta}(0)\cap \Rn=\{z\in\Rn\,;\:|z|=\delta \}$ and $\d 'B^+_{\delta}(0)=B^+_{\delta}(0)\cap \d\Rn=\{z\in\d\Rn\,;\:|z|<\delta \}$. Thus, $\d B^+_{\delta}(0)=\d 'B^+_{\delta}(0)\cup \d^+B^+_{\delta}(0)$.

In various parts of the text, we will make use of Fermi coordinates $$\psi:B^+_{\delta}(0)\to M$$ centered at a point $x_0\in\d M$. In this case, we will work in $B^+_{\delta}(0)\subset \Rn$. 

We will denote by $D_{\delta}(x_0)$ the metric ball on $\d M$ (with respect to the induced metric $\bar{g}$) of radius $\delta>0$ centered at $x_0\in\d M$.
The volume forms of $M$ and $\d M$ will be denoted by $dv_g$ and $d\sigma_g$, respectively. By $\eta$ we will denote the inward unit normal vector to $\d M$.
The n-dimensional sphere of radius $r$ in $\R^{n+1}$ will be denoted by $S_r^n$.
Finally, $\sigma_{n}$ will denote the area of the n-dimensional unit sphere $S^n_1$.  

\subsection{Standard solutions in the Euclidean half-space}\label{subsec:standard_solutions}

In this section we will study the Euclidean Yamabe equation in $\Rn$ and its linearization.

The simplest example of solution to the Yamabe-type problem we are concerned is the ball in $\R^n$ with the canonical Euclidean metric. This ball is conformally equivalent to the half-space $\Rn$ by the inversion 
$F:\mathbb{R}_+^n\to B^n\backslash\{ (0,...,0,-1)\}$
with respect to the sphere $S^{n-1}_1(0,...,0,-1)$ with center $(0,...,0,-1)$ and radius $1$. Here, $B^n=B_{1/2}(0,...,0,-1/2)$ is the Euclidean ball in $\R^{n}$ with center $(0,...,0,-1/2)$ and radius $1/2$. The expression for $F$ is
$$F(y_1,...y_n)=\frac{(y_1,...,y_{n-1},y_n+1)}{y_1^2+...+y_{n-1}^2+(y_n+1)^2}+(0,...,0,-1)\,,$$
and of course  its inverse mapping $F^{-1}$ has the same expression.
An easy calculation shows that $F$ is a conformal map and $F^*g_{eucl} = U^{\frac{4}{n-2}}g_{eucl}$ in $\Rn$, where $g_{eucl}$ is the Euclidean metric and  
$U(y)=(y_1^2+...+y_{n-1}^2+(y_n+1)^2)^{-\frac{n-2}{2}}$.
The function $U$ satisfies
\begin{equation}
\label{eq:U}
\begin{cases}
\Delta U = 0\,,&\text{in}\:\mathbb{R}_+^n\,,\\
\frac{\partial U}{\partial y_n}+(n-2)U^{\frac{n}{n-2}}=0\,,&\text{on}\:\partial\mathbb{R}_+^n\,.
\end{cases}
\end{equation}

Since the equations (\ref{eq:U}) are invariant by horizontal translations and scalings with respect to the origin, we obtain the following family of solutions to the equation (\ref{eq:U}):
\begin{equation}\label{fam:U}
U_{\l,z}(y)=\left(\frac{\l}{(\l+y_n)^2+\sum_{j=1}^{n-1}(y_j-z_j)^2}\right)^{\frac{n-2}{2}}\,,
\end{equation}
where $\l>0$ and $z=(z_1,...,z_{n-1})\in \R^{n-1}$.

In fact, the converse statement is also true: by a Liouville-type theorem in \cite{li-zhu} (see also \cite{escobar1} and \cite{chipot-shafrir-fila}), any non-negative solution to the equations (\ref{eq:U}) is of the form (\ref{fam:U}) or is identically zero. 


The existence of the family of solutions (\ref{fam:U}) has two important consequences. First, we see that the set of solutions to the equations (\ref{eq:U}) is non-compact. In particular, the set of solutions to the equations (\ref{eq:u'}) is not compact when $M^n$ is conformally equivalent to $B^n$. Secondly, the functions $\frac{\d U}{\d y_j}$, for $j=1,...,n-1$, and $\frac{n-2}{2}U+y^b\frac{\d U}{\d y^b}$, are solutions to the following homogeneous linear problem:
\ba
\begin{cases}\label{linear:homog}
\Delta\psi=0\,,&\text{in}\:\Rn\,,
\\
\frac{\d\psi}{\d y_n}+nU^{\frac{2}{n-2}}\psi=0\,,&\text{on}\:\d\Rn\,.
\end{cases}
\end{align}
\begin{notation}
We set $J_j=\frac{\d U}{\d y_j}$, for $j=1,...,n-1$, and $J_n=\frac{n-2}{2}U+y^b\frac{\d U}{\d y^b}$.
\end{notation}

Now, we will show that linear combinations of $J_1,...,J_n$ are the only solutions to the equations (\ref{linear:homog}) under a certain decay hypothesis. This result is similar to the one obtained in \cite{chen-lin} for the case of manifols without boundary. More precisely we have:
\begin{lemma}
\label{classifLinear}
Suppose $\psi$ is a solution to
\begin{equation}
\label{eq:Ulinear}
\begin{cases}
\Delta \psi = 0\,,&\text{in}\:\mathbb{R}_+^n\,\\
\frac{\partial \psi}{\partial y_n}+nU^{\frac{2}{n-2}}\psi=0\,,&\text{on}\:\partial \mathbb{R}_+^n\,.
\end{cases}
\end{equation}
If $\psi(y)=O((1+|y|)^{-\alpha})$ for some $\alpha>0$, then there exist constants $c_1,...,c_{n}$ such that
\[
\psi (y)=\sum_{a=1}^{n} c_a J_a(y)\,.
\]
\end{lemma}

The following result will be used in the proof of Lemma \ref{classifLinear}:
\begin{lemma}\label{eigenvalues:Bn}
The eigenvalues $\lambda$ of the problem
\begin{equation}\label{linear:3'}
\begin{cases}
\Delta \bar{\psi} = 0\,,&\text{in}\:B^n\,,\\
\frac{\partial \bar{\psi}}{\partial \eta}+\lambda\bar{\psi}=0\,,&\text{on}\:\partial B^n
\end{cases}
\end{equation}
are given by $\{\lambda_k=2k\}_{k=0}^{\infty}$. (Recall that $\eta$ points inwards.) The corresponding eigenvectors are the harmonic homogeneous polynomials of degree $k$ restricted to $B^n$. Here, the coefficients of the polynomials are given by the coordinate functions of $\R^n$ with center $(0,...,0,-1/2)$. In particular, the constant function $1$ generates the eigenspace associated to the eigenvalue $\l_0=0$ and the coordinate functions $z_1,...,z_n$ restricted to $B^n$ generate the eigenspace associated to the eigenvalue $\l_1=2$. 

Moreover, $z_j\circ F=\frac{-1}{n-2}U^{-1}J_j$, for $j=1,...,n-1$, and $z_n\circ F=\frac{1}{n-2}U^{-1}J_n$.
\end{lemma}

\bp
The first part is an easy consequence of the fact that the spherical harmonics generate $L^2(S^{n-1})$. The last part is a straightforward computation. 
\ep

\begin{proof}[Lemma \ref{classifLinear}]
Recall that the conformal Laplacian satisfies
\begin{equation}\label{propr:L}
L_{\zeta^{\frac{4}{n-2}}g}(\zeta^{-1}u)=\zeta^{-\frac{n+2}{n-2}}L_{g}u,
\end{equation}
for any smooth functions $\zeta>0$ and $u$. Similarly, the boundary operator $B_g$ satisfies
\begin{equation}\label{propr:B}
B_{\zeta^{\frac{4}{n-2}}g}(\zeta^{-1}u)=\zeta^{-\frac{n}{n-2}}B_{g}u\,.
\end{equation}
Hence, the equations \eqref{eq:Ulinear} are equivalent to
\begin{equation}\notag
\begin{cases}
\Delta \bar{\psi} = 0\,,&\text{in}\:B^n\backslash\{ (0,...,0,-1)\}\,,
\\
\frac{\partial \bar{\psi}}{\partial \eta}+2\bar{\psi}=0\,,&\text{on}\:\partial B^n\backslash\{ (0,...,0,-1)\}\,,
\end{cases}
\end{equation}
where $\bar{\psi}=(U^{-1}\psi)\circ F^{-1}$.
The hypothesis $\psi(y)=O((1+|y|)^{-\alpha})$, $0<\alpha <n-2$ implies that $\bar{\psi}\in L^p(B^n)$, for any $\frac{n}{n-2}<p<\frac{n}{n-2-\alpha}$. Lemma \ref{extensao:sol} ensures that $\bar{\psi}$ is a weak solution to 
\begin{equation}\notag
\begin{cases}
\Delta \bar{\psi} = 0\,,&\text{in}\:B^n\,,
\\
\frac{\partial \bar{\psi}}{\partial \eta}+2\bar{\psi}=0\,,&\text{on}\:\partial B^n\,.
\end{cases}
\end{equation}
It follows from elliptic theory that $\bar{\psi}\in C^{\infty}(B^n)$. In other words, $\psi$ is a solution to the equations \eqref{eq:Ulinear} if and only if $\bar{\psi}$ is an eigenfunction associated to the first nontrivial eigenvalue $\lambda_1=2$ of the problem (\ref{linear:3'}). The result now follows from Lemma \ref{eigenvalues:Bn}.
\end{proof}


\subsection{Coordinate expansions for the metric}\label{subsec:coodinate}

In this section we will write expansions for the metric $g$ in Fermi coordinates. We will also recall the conformal Fermi coordinates, introduced by Marques in \cite{coda1}, that will simplify the computations in the next sections. The conformal Fermi coordinates play the same role that the conformal normal coordinates (see \cite{lee-parker}) did in the case of manifolds without boundary.
\begin{definition}
Let $x_0\in\d M$. We choose boundary geodesic normal coordinates $(z_1,...,z_{n-1})$, centered at $x_0$, of the point $x\in\d M$.
We say that $z=(z_1,...,z_n)$, for small $z_n\geq 0$, are the {\it{Fermi coordinates}} (centered at $x_0$) of the point $\exp_{x}(z_n\eta(x))\in M$. Here, we denote by $\eta(x)$ the inward unit normal vector to $\d M$ at $x$. In this case, we have a map defined by $\psi(z)=\exp_{x}(z_n\eta(x))$, taking values in a subset of $\Rn$.
\end{definition}

It is easy to see that in these coordinates $g_{nn}\equiv 1$ and $g_{jn}\equiv 0$, for $j=1,...,n-1$.

We fix $x_0\in\d M$. We use Fermi coordinates $\psi:B^+_{\delta}(0)\to M$ centered at $x_0$ and work in $B_{\delta}^+(0)\subset\Rn$, for some small $\delta>0$.  
\begin{notation}
We set
$$
|\d^kg|=\max_{z\in B^+_{\delta}(0)}\sum_{|\a|= k}\sum_{a,b=1}^{n}|\d^{\a}g_{ab}|(\psi(z))\,,
$$
where $\a$ denotes a multiindex.
We write $|\d g|=|\d^1g|$ for short. 
\end{notation}

The following proposition gives the expansion for the Riemannian metric $g$ in Fermi coordinates:
\begin{proposition}\label{exp:g}
For $z=(z_1,...,z_n)\in B_{\delta}^+(0)$,
\begin{align}
g^{ij}(\psi(z))=\;&\delta_{ij}+2h_{ij}(x_0)z_n
+\frac{1}{3}\bar{R}_{ikjl}(x_0)z_kz_l+2h_{ij;\,\underline{k}}(x_0)z_nz_k\notag
\\
&+(R_{ninj}+3h_{ik}h_{kj})(x_0)z_n^2
+O(|\d^3g||z|^3).\notag
\end{align}
\end{proposition}

\bp
This is proved as in Lemma 2.2 of \cite{coda1}.
\ep

The existence of conformal Fermi coordinates and some of its consequences are stated as follows:
\begin{proposition}\label{conf:fermi:thm}
For any given integer $N\geq 1$, there is a metric $\tilde{g}$, conformal to $g$, such that in $\tilde{g}$-Fermi coordinates $\tilde{\psi}:B^+_{\delta}(0)\to M$ centered at $x_0$, we have
$$
(\det \tilde{g})(\tilde{\psi}(z))=1+O(|z|^N)\,.
$$  
Moreover, $\tilde{g}$ can be written as $\tilde{g}=fg$, where $f$ is a positive function with $f(x_0)=1$ and $\frac{\d f}{\d z_k}(x_0)=0$ for $k=1,...,n-1$. In this metric we also have
\\\\
(i) $\bar{R}_{ij}(x_0)=\bar{R}_{ij;\,\underline{k}}(x_0)=0$;
\\
(ii) $R_{nn}(x_0)+(h_{ij})^2(x_0)=0$;
\\
(iii) $h(\tilde{\psi}(z))=O(|z|^N)$, where $N$ can be taken arbitrarily large.
\end{proposition}

\bp
The first part is Proposition 3.1 of \cite{coda1}. Items (i) and (ii) are proved as in Proposition 3.2 of the same paper. Item (iii) follows from the fact that
\begin{equation*}
h_g=\frac{-1}{2(n-1)}g^{ij}g_{ij,\,n}=\frac{-1}{2(n-1)}(\log\det g)_{,\,n}\,.
\end{equation*}
\ep


\subsection{Conformal scalar and mean curvature equations}\label{subsec:conformal_scalar}
In this section we will introduce the partial differential equation we will work with in the next sections. We will also discuss some of its properties related to conformal deformation of metrics.

Let $u$ be a positive smooth solution to
\begin{align}\label{eq:u:geral}
\begin{cases}
L_{g}u=0,&\text{in}\:M\,,
\\
B_{g}u+(n-2)f^{-\tau}u^{p}=0,&\text{on}\:\d M\,,
\end{cases}
\end{align}
where $\tau=\frac{n}{n-2}-p$, $1+\gamma_0\leq p\leq\frac{n}{n-2}$ for some fixed $\gamma_0>0$ and $f$ is a positive function. 
\begin{notation}
We say that $u\in\mathcal{M}_p$ if $u$ is a positive smooth solution to the equations (\ref{eq:u:geral}).
\end{notation}

The equations (\ref{eq:u:geral}) have an important scaling invariance property. We fix $x_0\in \d M$ and let $\delta>0$ be small. We consider Fermi coordinates $\psi:B^+_{\delta}(0)\to M$ centered at $x_0$. Given $s>0$ we define the renormalized function 
$$
v(y)=s^{\frac{1}{p-1}}u(\psi(sy))\,,\:\:\:\text{for}\: y\in B^+_{\delta s^{-1}}(0)\,.
$$ 
Then
\begin{align}\notag
\begin{cases}
L_{\hat{g}}v=0,&\text{in}\:B^+_{\delta s^{-1}}(0)\,,
\\
B_{\hat{g}}v+(n-2)\hat{f}^{-\tau}v^{p}=0,&\text{on}\: \d 'B^+_{\delta s^{-1}}(0)\,,
\end{cases}
\end{align}
where $\hat{f}(y)=f(\psi(sy))$ and the coefficients of the metric $\hat{g}$ in Fermi coordinates are given by $\hat{g}_{kl}(y)=g_{kl}(\psi(sy))$.

The reason to work with the equations (\ref{eq:u:geral}) instead of the equations (\ref{main:eq}) is that the first one has an important conformal invariance property.  Suppose $\tilde{g}=\zeta^{\frac{4}{n-2}}g$ is a metric conformal to $g$.
It follows from the properties \eqref{propr:L} and \eqref{propr:B} that, if $u$ is a solution to the equations (\ref{eq:u:geral}), then $\zeta^{-1}u$ satisfies
\begin{align}
\begin{cases}\notag
L_{\tilde{g}}(\zeta^{-1}u)=0,&\text{in}\:M\,,
\\
B_{\tilde{g}}(\zeta^{-1}u)+(n-2)(\zeta f)^{-\tau}(\zeta^{-1}u)^{p}=0,&\text{on}\:\d M\,,
\end{cases}
\end{align}
which is again equations of the same type. 
\begin{notation}
Let $\Omega\subset M$ be a domain in a Riemannian manifold $(M,g)$.
Let $\{g_i\}$ be a sequence of metrics on $\Omega$. We say that $u_i\in\mathcal{M}_i$ if $u_i>0$ satisfies
\begin{align}\label{eq:ui}
\begin{cases}
L_{g_i}u_i=0,&\text{in}\:\Omega\,,
\\
B_{g_i}u_i+(n-2)f_i^{-\tau_i}u_i^{p_i}=0,&\text{on}\:\Omega\cap \d M\,,
\end{cases}
\end{align}
where $\tau_i=\frac{n}{n-2}-p_i$ and $1+\gamma_0\leq p_i\leq\frac{n}{n-2}$ for some fixed $\gamma_0>0$. 
\end{notation}

In many parts of this article we will work with sequences $\{u_i\in\mathcal{M}_i\}_{i=1}^{\infty}$. In this case, we assume that $f_i\to f$ in the $C^1_{loc}$ topology, for some positive function $f$, and that $g_i\to g_0$ in the $C^3_{loc}$ topology, for some metric $g_0$. 

By the conformal invariance stated above, we are allowed to replace the metric $g_i$ by $\zeta_i^{\frac{4}{n-2}}g_i$ as long as we have control of the conformal factors $\zeta_i$. In this case, we replace the sequence $\{u_i\}$ by $\{\zeta_i^{-1}u_i\}$. In particular, we can use conformal Fermi coordinates centered at some point $x_i\in\d M$.


\section{A Pohozaev-type identity}

In this section we prove the Pohozaev-type identity we will use in the subsequent blow-up analysis.

Recall that we have denoted $B^+_{\delta}=B^+_{\delta}(0)\subset\Rn$ for short (see Section \ref{subsec:notations}).


\begin{proposition}\label{Pohozaev}

Let $u$ be a solution to
\begin{equation}\notag
\begin{cases}
\Delta_{g}u-\frac{n-2}{4(n-1)}R_gu=0\,,&\text{in}\:B^+_{\delta}\,,\\
\partial_n u-\frac{n-2}{2} h_gu+Kf^{-\tau}u^{p}=0\,,&\text{on}\:\partial' B^+_{\delta} \,,
\end{cases}
\end{equation}
where $K$ is a constant and $g$ is a metric on $B^+_{\delta}$.
Let $0<r<\delta$. We set
$$P(u,r)=\int_{\partial^+ B^+_r}\left(\frac{n-2}{2}u\frac{\partial u}{\partial r}-\frac{r}{2}|\nabla u|^2+r\left|\frac{\partial u}{\partial r}\right|^2\right)d\sigma_r+\frac{r}{p+1}\int_{\partial\,(\partial' B^+_{r})}Kf^{-\tau}u^{p+1}d\bar{\sigma}_r\,.$$
Then
\begin{align}
P(u,r)&=-\int_{B_r^+}\left(z^a\partial_a u+\frac{n-2}{2}u\right)A_g(u)dz
+\frac{n-2}{2}\int_{\partial' B^+_{r}}\left(\bar{z}^k\partial_k u+\frac{n-2}{2}u\right)h_gud\bar{z}\notag
\\
&-\frac{\tau}{p+1}\int_{\partial' B^+_{r}}K(\bar{z}^k\partial_k f)f^{-\tau-1}u^{p+1}d\bar{z}
+\left(\frac{n-1}{p+1}-\frac{n-2}{2}\right)\int_{\partial' 
B^+_{r}}Kf^{-\tau}u^{p+1}d\bar{z}\,,\notag
\end{align}
where $A_g=\Delta_g-\Delta-\frac{n-2}{4(n-1)}R_g$. Here, $\Delta$ stands for the Euclidean Laplacian and $\nabla$ for the Euclidean gradient.
\end{proposition}

\begin{proof}
Observe that, for each $a=1,...,n$ fixed, integrating by parts we have
\begin{align}
\int_{B^+_r}(z^b\partial_bu)\d_{a}\d_{a}udz
&+\int_{B^+_r}\delta^{ab}(\partial_bu)(\partial_au)dz
+\frac{1}{2}\int_{B^+_r}z^b\partial_b(\partial_au)^2dz\notag
\\
&=\frac{1}{r}\int_{\partial^+B^+_r}(z^b\partial_bu)(z^a\partial_au)d\sigma_r
-\int_{\partial 'B^+_r}(\bar{z}^k\partial_k u)(\partial_au)\delta^{a}_nd\bar{z}\,.\notag
\end{align}
Summing in $a=1,...,n$ we obtain
\begin{align}\label{poho1}
\int_{B^+_r}(z^b\partial_bu)\Delta udz
+\int_{B^+_r}|\nabla u|^2dz
&+\frac{1}{2}\sum_{a}\int_{B^+_r}z^b\partial_b(\partial_a u)^2dz\notag
\\
&=r\int_{\partial^+B^+_r}\left|\frac{\partial u}{\partial r} \right|^2d\sigma_r
-\int_{\partial 'B^+_r}(\bar{z}^k\partial_k u)(\partial_n u)d\bar{z}\,.
\end{align}
On the other hand, integrating by parts, we have
\begin{align}\label{poho2}
\frac{1}{2}\sum_{a}\int_{B^+_r}z^b\partial_b(\partial_au)^2dz
&=-\frac{n}{2}\sum_{a}\int_{B^+_r}(\partial_a u)^2dz
+\frac{r}{2}\sum_{a}\int_{\partial^+B^+_r}(\partial_a u)^2d\sigma_r\notag
\\
&\hspace{1cm}-\frac{1}{2}\sum_{a}\int_{\partial 'B^+_r}z^b\delta_{b}^n(\partial_a u)^2d\bar{z}\notag
\\
&=-\frac{n}{2}\int_{B^+_r}|\nabla u|^2dz
+\frac{r}{2}\int_{\partial^+B^+_r}|\nabla u|^2d\sigma_r
\end{align}
and
\begin{align}\label{poho3}
\int_{\partial 'B^+_r}(\bar{z}^k\partial_k u)(\partial_n u)d\bar{z}
&=-\int_{\partial 'B^+_r}(\bar{z}^k\partial_k u)(Kf^{-\tau}u^{p}-\frac{n-2}{2} h_gu)d\bar{z}\notag
\\
&=-\frac{1}{p+1}\int_{\partial 'B^+_r}K\bar{z}^k\partial_k(u^{p+1})f^{-\tau}d\bar{z}\notag
\\
&\hspace{1cm}+\frac{n-2}{2}\int_{\partial 'B^+_r}(\bar{z}^k\partial_k u)h_gud\bar{z}\notag
\\
&=\frac{n-1}{p+1}\int_{\partial 'B^+_r}Kf^{-\tau}u^{p+1}d\bar{z}
+\frac{1}{p+1}\int_{\partial 'B^+_r}K(\bar{z}^k\partial_k f^{-\tau})u^{p+1}d\bar{z}\notag
\\
&\hspace{1cm}-\frac{r}{p+1}\int_{\partial (\partial 'B^+_r)} Kf^{-\tau}u^{p+1}d\bar{\sigma}_r\notag
\\
&\hspace{2cm}+\frac{n-2}{2}\int_{\partial 'B^+_r}(\bar{z}^k\partial_k u)h_gud\bar{z}\,.
\end{align}
Substituting equalities \eqref{poho2} and \eqref{poho3} in \eqref{poho1} we obtain
\begin{align}\label{poho4}
\int_{B^+_r}(z^b\partial_b u)&\Delta u dz
-\frac{n-2}{2}\int_{B^+_r}|\nabla u|^2dz
+\frac{r}{2}\int_{\partial^+B^+_r}|\nabla u|^2d\sigma_r\notag
\\
&=r\int_{\partial^+B^+_r}\left|\frac{\partial u}{\partial r} \right|^2d\sigma_r
-\frac{n-1}{p+1}\int_{\partial 'B^+_r}Kf^{-\tau}u^{p+1}d\bar{z}\notag
\\
&\hspace{1cm}-\frac{1}{p+1}\int_{\partial 'B^+_r}K(\bar{z}^k\partial_k f^{-\tau})u^{p+1}d\bar{z}
+\frac{r}{p+1}\int_{\partial (\partial 'B^+_r)}Kf^{-\tau}u^{p+1}d\bar{\sigma}_r\notag
\\
&\hspace{2cm}-\frac{n-2}{2}\int_{\partial 'B^+_r}(\bar{z}^k \partial_k u)h_gud\bar{z}\,.
\end{align}
Using
\begin{equation}\notag
\int_{B^+_r}|\nabla u|^2dz
=-\int_{B^+_r}u\Delta udz 
+\int_{\partial^+B^+_r}u\frac{\partial u}{\partial r}d\sigma_r
+\int_{\partial 'B^+_r}(Kf^{-\tau}u^{p+1}-\frac{n-2}{2} h_gu^2)d\bar{z}
\end{equation}
and $\Delta u=-A_g(u)$ in equality \eqref{poho4} we get the result.
\end{proof}



\section{Isolated and isolated simple blow-up points}\label{sec:isolated:simple}

In this section we will discuss the notions of isolated and isolated simple blow-up points and prove some of their properties. These notions are slight modifications of the ones used by  Felli and Ould Ahmedou in \cite{ahmedou-felli1} and \cite{ahmedou-felli2} and are inspired by similar definitions in the case of manifolds without boundary. 
 \begin{definition}\label{def:blow-up}
Let $\Omega\subset M$ be a domain in a Riemannian manifold $(M,g)$.
We say that $x_0\in\Omega\cap\d M$ is a {\it{blow-up point}} for the sequence $\{u_i\in\mathcal{M}_i\}_{i=1}^{\infty}$, if there is a sequence $\{x_i\}\subset\Omega\cap\d M$ such that 

(1) $x_i\to x_0$;

(2) $u_i(x_i)\to\infty$; 

(3) $x_i$ is a local maximum of $u_i|_{\d M}$.\\
Briefly we say that $x_i\to x_0$ is a blow-up point for $\{u_i\}$. The sequence $\{u_i\}$ is called a {\it{blow-up sequence}}.
\end{definition}

\noindent
{\bf{Convention}} If $x_i\to x_0$ is a blow-up point, we use $g_i$-Fermi coordinates $$\fermi$$ centered at $x_i$ (see Section \ref{subsec:coodinate}) and work in $B^+_{\delta}(0)\subset\Rn$, for some small $\delta>0$.
\begin{notation}
If $x_i\to x_0$ is a blow-up point we set $M_i=u_i(x_i)$, $\ei=M_i^{-(p_i-1)}$.
\end{notation}


\subsection{Isolated blow-up points}

We define the notion of an isolated blow-up point as follows:
\begin{definition}\label{def:isolado}
We say that a blow-up point $x_i\to x_0$ is an {\it{isolated}} blow-up point for $\{u_i\}$ if there exist $\delta,C>0$ such that 
\begin{equation}\notag
u_i(x)\leq Cd_{\bar{g}_i}(x,x_i)^{-\frac{1}{p_i-1}}\,,\:\:\:\:\text{for all}\: x\in\d M\backslash \{x_i\}\,,\:d_{\bar{g}_i}(x,x_i)< \delta\,.
\end{equation}
(We recall that $\bar{g}_i$ denotes the induced metric on the boundary.)
\end{definition}
\begin{remark}
Since Fermi coordinates are normal on the boundary, the above definition is equivalent to
\begin{equation}\label{des:isolado}
u_i(\psi_i(z))\leq C|z|^{-\frac{1}{p_i-1}}\,,\:\:\:\:\text{for all}\: z\in \d 'B^+_{\delta}(0)\backslash\{0\}\,.
\end{equation}
\end{remark}
\begin{remark}\label{scaling:inv}
Note that the definition of an isolated blow-up point is invariant under renormalization, which was descrided in Section 2.4. This follows from the fact that if $v_i(y)=s^{\frac{1}{p_i-1}}u_i(\psi_i(sy))$, then
$$
u_i(\psi_i(z))\leq C|z|^{-\frac{1}{p_i-1}}\Longleftrightarrow v_i(y)\leq C|y|^{-\frac{1}{p_i-1}}\,,
$$
where $z=sy$.
\end{remark}


The first result concerning isolated blow-up points states that the inequality (\ref{des:isolado}) also holds for points $z\in B^+_{\delta}(0)\backslash\{0\}$.
\begin{lemma}\label{Harnack:bordo}
Let $x_i\to x_0$ be an isolated blow-up point. Then $\{u_i\}$ satisfies
\begin{equation}\notag
u_i(\psi_i(z))\leq C|z|^{-\frac{1}{p_i-1}}\,,\:\:\:\text{for all}\: z\in B^+_{\delta}(0)\backslash\{0\}\,.
\end{equation}
\end{lemma}
\begin{proof} Let $0<s<\frac{\delta}{3}$ and set $v_{i}(y)=s^{\frac{1}{p_i-1}}u_i(\psi_i(sy))$ for $|y|< 3$.
Then $v_i$ satisfies
\begin{align}
\begin{cases}
L_{\hat{g}_i}v_i=0,&\text{in}\:B^+_3(0),\notag
\\
(B_{\hat{g}_i}+(n-2)\hat{f}_i^{-\tau_i}v_i^{p_i-1})v_i=0,&\text{on}\:\partial 'B^+_3(0)\,,\notag
\end{cases}
\end{align}
where $(\hat{g}_i)_{kl}(y)=(g_i)_{kl}(\psi_i(sy))$ and $\hat{f}_i(y)=f_i(\psi_i(sy))$.
By the scaling invariance (Remark \ref{scaling:inv}) $v_i$ is uniformly bounded in compact subsets of $\d 'B^+_3(0)\backslash\{0\}$. Hence, Lemma \ref{Harnack:han-li} and interior Harnack estimates give
\begin{equation}\label{Harnack:bordo:1}
\max_{B_{2}^+(0)\backslash B_{1/2}^+(0)} v_i\leq C\min_{B_{2}^+(0)\backslash B_{1/2}^+(0)} v_i\,.
\end{equation}
The result now follows from the inequality (\ref{Harnack:bordo:1}).
\end{proof}

A corollary of the proof of Lemma \ref{Harnack:bordo} is the following Harnack-type  inequality:
\begin{lemma}\label{Harnack}
Let $x_i\to x_0$ be an isolated blow-up point and $\delta$ as in Definition \ref{def:isolado}. 
Then there exists $C>0$ such that for any $0<s<\frac{\delta}{3}$ we have
\begin{equation}\notag
\max_{B_{2s}^+(0)\backslash B_{s/2}^+(0)} (u_i\circ\psi_i)\leq C\min_{B_{2s}^+(0)\backslash B_{s/2}^+(0)} (u_i\circ\psi_i)\,.
\end{equation}
\end{lemma}



The  next proposition says that, in the case of an isolated blow-up point, the sequence $\{u_i\}$, when renormalized, converges to the standard Euclidean solution $U$ (see Section \ref{subsec:standard_solutions}).
 \begin{proposition}\label{form:bolha}
Let $x_i\to x_0$ be an isolated blow-up point. We set
$$
v_i(y)=M_i^{-1}(u_i\circ\psi_i)(M_i^{-(p_i-1)} y)\,,\:\:\:\: \text{for}\: y\in B^+_{\delta M_i^{p_i-1}}(0)\,.
$$
Then given $R_i\to\infty$ and $\b_i\to 0$, after choosing subsequences, we have  
\\\\
(a) $|v_i-U|_{C^2(B^+_{R_i}(0))}<\b_i$;
\\
(b) $\lim_{i\to\infty}\frac{R_i}{\log M_i}= 0$;
\\
(c) $\lim_{i\to\infty}p_i=\frac{n}{n-2}$.
\end{proposition}

The proof of Proposition \ref{form:bolha} is analogous to Lemma 2.6 of \cite{ahmedou-felli1} or Proposition 4.3 of \cite{marques}. It uses the fact that, by the Liouville-type theorems of \cite{hu} and \cite{li-zhu}, every non-negative solution to 
\ba\label{eq:v}
\begin{cases}
\Delta v=0,&\text{in}\:\Rn,
\\
\d_n v+(n-2)v^{p_0}=0,&\text{on}\:\d\Rn,
\end{cases}
\end{align}
for $1<p_0\leq\frac{n}{n-2}$, is either identically zero or is 
of the form (\ref{fam:U}), in which case $p_0=\frac{n}{n-2}$.
\begin{remark}\label{rk:mud_conf}
Let $x_i\to x_0$ and consider a conformal change $\zeta_i^{\frac{4}{n-2}}g_i$ of the metrics $g_i$ (see the last paragraph of Section \ref{subsec:conformal_scalar}). Suppose that the conformal factors $\zeta_i>0$ are uniformly bounded (above and below) with $\zeta_i(x_i)=1$ and $\frac{\d \zeta_i}{\d z_k}(x_i)=0$ for $k=1,...,n-1$.
Then, once we have proved Proposition \ref{form:bolha}, it is not difficult to see that $x_i\to x_0$ is an isolated blow-up point for $\{u_i\}$ if and only it is for $\{\zeta_i^{-1}u_i\}$. This is the case when we use conformal Fermi coordinates (see Proposition \ref{conf:fermi:thm}) centered at $x_i$.
\end{remark}



The following lemma will be used later when we consider the  set of blow-up points.
\begin{lemma}\label{isolado:compactos}
Given $R,\beta>0$,  there exists $C_0=C_0(R,\beta)>0$ such that if $u\in\mathcal{M}_p$ and $S\subset\d M$ is a compact set, we have the following:

If $\max_{x\in\d M\backslash S}\left(u(x)d_{\bar{g}}(x,S)^{\frac{1}{p-1}}\right)\geq C_0$, then $\frac{n}{n-2}-p<\beta$ and there exists $x_0\in \d M\backslash S$, local maximum of $u$, such that 
\begin{equation}\label{isolado:compactos:0}
\left|u(x_0)^{-1}u(\psi(z))-U(u(x_0)^{p-1}z)\right|_{C^2(B^+_{2r_0}(0))}<\b\,,
\end{equation}
where $r_0=Ru(x_0)^{-(p-1)}$.
Here, we are using Fermi coordinates $\psi:B_{2r_0}^+(0)\to M$ centered at $x_0\in\partial M$. 
If $\emptyset$ is the empty set, we define $d_{\bar{g}}(x,\emptyset)=1$. 
\end{lemma}
\begin{proof}
Suppose by contradiction there exist $R,\b >0$ satisfying the following: for all $C_0>0$, there exist $p\in\left(1,\frac{n}{n-2}\right]$, $u\in\mathcal{M}_p$ and a compact set $S\subset\d M$ such that 
$$\max_{x\in \d M\backslash S}\left(u(x)d_{\bar{g}}(x,S)^{\frac{1}{p-1}}\right)\geq C_0$$ 
holds and either $p\leq\frac{n}{n-2}-\beta$ or no such point $x_0$ exists. Hence, we can suppose that there are sequences $p_i\in\left(1,\frac{n}{n-2}\right]$, $u_i\in\mathcal{M}_{p_i}$ and
$$w_i(x'_i)=\max_{x\in\d M\backslash\ S_i} w_i(x)\to \infty\,,\:\:\:\:\text{where}\:\: w_i(x)=u_i(x)d_{\bar{g}}(x,S_i)^{\frac{1}{p_i-1}}\,.$$  
Here, $x'_i\in\d M$ and $S_i$ is a compact subset of $\d M$. We assume that $p_i\to p_0$, for some $p_0\in\left(1,\frac{n}{n-2}\right]$, and $x_i'\to x_0'$ for some $x'_0\in\d M$. We set $N_i~=~u_i(x_i')$. Observe that $N_i\to\infty$.
 
We use Fermi coordinates $\fermi$ centered at $x'_i$ and set 
$$
v_i(y)=N_i^{-1}(u_i\circ\psi_i)(N_i^{-(p_i-1)}y)\,,\:\:\:\:\text{for}\:y\in B^+_{\delta N_i^{p_i-1}}(0)\,. 
$$
It follows from the discussion in Section 2.4 that $v_i$ satisfies
\begin{equation}\notag
\begin{cases}
L_{\hat{g}_i}v_i=0\,,&\text{in}\:B^+_{\delta N_i^{p_i-1}}(0)\,,
\\
B_{\hat{g}_i}v_i+(n-2)\hat{f}_i^{-\tau_i}v_i^{p_i}=0\,,&\text{on}\:\d 'B^+_{\delta N_i^{p_i-1}}(0)\,,
\end{cases}
\end{equation}
where $\hat{f}_i(y)=f(\psi_i(N_i^{-(p_i-1)} y))$ and $\hat{g}_i$ stands for the metric with coefficients $(\hat{g}_i)_{kl}(y)=g_{kl}(\psi_i(N_i^{-(p_i-1)}y))$.
\\\\
{\it{Claim}} $v_i\leq C$ in compact subsets of $\Rn$.

\vspace{0.2cm}
Let $z\in\d 'B_{\delta}^+(0)$. Since $w_i(\psi_i(z))\leq w_i(x'_i)$, we have
$$\frac{d_{\bar{g}}(S_i,x'_i)-d_{\bar{g}}(x'_i,\psi_i(z))}{d_{\bar{g}}(S_i,x'_i)}
\leq \frac{d_{\bar{g}}(S_i,\psi_i(z))}{d_{\bar{g}}(S_i,x'_i)}
\leq\left(N_iu_i(\psi_i(z))^{-1}\right)^{p_i-1}\,.$$
On the other hand, 
$$\frac{d_{\bar{g}}(S_i,x'_i)-d_{\bar{g}}(x'_i,\psi_i(z))}{d_{\bar{g}}(S_i,x'_i)}
=1-\frac{N_i^{-(p_i-1)}|y|}{d_{\bar{g}}(S_i,x'_i)}
=1-w_i(x'_i)^{-(p_i-1)}|y|=1-o_i(1)|y|,$$
where we have set $y=N_i^{p_i-1}z$. This proves that $v_i\leq C$ in compact subsets of $\d\Rn$. Now the Claim follows from Lemma \ref{Harnack:han-li}.

Hence, we can suppose that $v_i\to v$ in $C^2_{loc}(\Rn)$ for some $v>0$ satisfying 
\ba\notag
\begin{cases}
\Delta v=0,&\text{in}\:\Rn,
\\
\d_n v+(n-2)f(x'_0)^{p_0-\frac{n}{n-2}}\,v^{p_0}=0,&\text{on}\:\d\Rn
\end{cases}
\end{align}
and $v(0)=1$. It follows from the Liouville-type theorem of \cite{hu} that $p_0=\frac{n}{n-2}$. Hence, $v$ satisfies the equations (\ref{eq:v}) and, by the results in \cite{li-zhu}, it is of the form (\ref{fam:U}).
Hence, we can find $y_{(i)}\in\d 'B^+_{\delta N_i^{p_i-1}}(0)$ local maxima of $v_i$, such that $y_{(i)}\to (z_1,...,z_{n-1},0)\in\d\Rn$. Then $u_i$ satisfies the  estimate (\ref{isolado:compactos:0}), for $i$ large, with $x_0=\psi_i(N_i^{-(p_i-1)}y_{(i)})$. Since $N_i^{p_i-1}d_{\bar{g}_i}(x'_i,S_i)=w_i(x'_i)^{p_i-1}\to\infty$, we see that $x_0\notin S_i$ for $i$ large.
This is a contradiction.
\end{proof}


Once we have proved Lemma \ref{isolado:compactos}, the proof of the following proposition is analogous to  Proposition 5.1 of \cite{li-zhu2} (see also Lemma 3.1 of \cite{schoen-zhang} or Proposition 1.1 of \cite{han-li}):
\begin{proposition}\label{conj:isolados}
Given small $\b>0$ and large $R>0$ there exist constants $C_0, C_1>0$, depending only on $\b$, $R$ and $(M^n,g)$, such that if $u\in\mathcal{M}_p$ and $\max_{\d M} u\geq C_0$, then $\frac{n}{n-2}-p<\b$ and
there exist $x_1,...,x_N\in\d M$, $N=N(u)\geq 1$, local maxima of $u$, such that:
\\\\
(1) If $r_j=Ru(x_j)^{-(p-1)}$ for $j=1,...,N$,  then $\{D_{r_j}(x_j)\subset\d M\}_{j=1}^{N}$ is a disjoint collection.
(We recall that $D_{r_j}(x_j)$ is the boundary metric ball (see Section \ref{subsec:notations}).)
\\\\
(2) For each $j=1,...,N$, 
$\:\:\:\:
\left|u(x_j)^{-1}u(\bar{\psi}_j(z))-U(u(x_j)^{p-1}z)\right|_{C^2(B^+_{2r_j}(0))}<\b
$,
\\
where we are using Fermi coordinates $\bar{\psi}_j:B^+_{2r_j}(0)\to M$ centered at $x_j$.
\\\\
(3) We have
$$
u(x)\,d_{\bar{g}}(x,\{x_1,...,x_N\})^{\frac{1}{p-1}}\leq C_1\,,\:\:\:\text{for all}\: x\in\d M\,,
$$
$$
u(x_j)\,d_{\bar{g}}(x_j,x_k)^{\frac{1}{p-1}}\geq C_0\,,\:\:\:\:\text{for any}\: j\neq k\,,\:j,k=1,...,N\,.
$$ 
\end{proposition}



\subsection{Isolated simple blow-up points}

Let us introduce the notion of an isolated simple blow-up point. Let $x_i\to x_0$ be an isolated blow-up point for $\{u_i\}$. Recall that we are using Fermi coordinates $\fermi$ centered at $x_i$.
We set 
$$\bar{u}_i(r)=\frac{2}{\sigma_{n-1}r^{n-1}}\int_{\partial^+B_r^+(0)}(u_i\circ\psi_i)d\sigma_r$$ 
and 
$w_i(r)=r^{\frac{1}{p_i-1}}\bar{u}_i(r)$, for $0<r<\delta$.

Note that the definition of $w_i$ is invariant under renormalization, which was descrided in Section 2.4. More precisely, if $v_i(y)=s^{\frac{1}{p_i-1}}u_i(\psi_i(sy))$, then
$$
r^{\frac{1}{p_i-1}}\bar{v}_i(r)=(sr)^{\frac{1}{p_i-1}}\bar{u}_i(sr)\,.
$$
\begin{definition}\label{def:simples}
An isolated blow-up point $x_i\to x_0$ for $\{u_i\}$ is {\it{simple}} if there exists $\delta>0$ such that $w_i$ has exactly one critical point in the interval $(0,\delta)$.
\end{definition}
\begin{remark}\label{rk:def:equiv} Let $x_i\to x_0$ be an isolated blow-up point and $R_i\to\infty$. Using Proposition \ref{form:bolha} it is not difficult to see that, choosing a subsequence, $r\mapsto r^{\frac{1}{p_i-1}}\bar{u}_i(r)$ has exactly one critical point in the interval $(0,r_i)$, where $r_i=R_iM_i^{-(p_i-1)}\to 0$. Moreover, its derivative is negative right after the critical point. Hence, if $x_i\to x_0$ is isolated simple then there exists $\delta>0$ such that $w_i'(r)<0$ for all $r\in [r_i,\delta)$. 
\end{remark}
\begin{notation}
In this section we define $$v_i(y)=M_i^{-1}(u_i\circ\psi_i)(M_i^{-(p_i-1)}y)\,,\:\:\:\:\:\text{for}\: y\in B_{M_i^{p_i-1}\delta}^+(0)\,.$$ 
\end{notation}
The next proposition is an important property of isolated simple blow-up points. 
\begin{proposition}\label{estim:simples}
Let $x_i\to x_0$ be an isolated simple blow-up point for $\{u_i\}$.  Then there exist $C,\delta>0$ such that
\\\\
(a) $M_iu_i(\psi_i(z))\leq C|z|^{2-n}$\:\:\: for all $z\in B^+_{\delta}(0)\backslash\{0\}$;
\\\\
(b) $M_iu_i(\psi_i(z))\geq C^{-1}G_i(z)$\:\:\: for all $z\in B^+_{\delta}(0)\backslash B^+_{r_i}(0)$, where $G_i$ is the Green's function so that:
\begin{align}
\begin{cases}
L_{g_i}G_i=0, &\text{in}\; B_{\delta}^+(0)\backslash\{0\},\notag
\\
G_i=0, &\text{on}\; \partial^+B_{\delta}^+(0),\notag
\\
B_{g_i}G_i=0, &\text{on}\;\partial 'B_{\delta}^+(0)\backslash\{0\}\notag 
\end{cases}
\end{align}
and $|z|^{n-2}G_i(z)\to 1$, as $|z|\to 0$. Here, $r_i$ is defined as in Remark \ref{rk:def:equiv}.
\end{proposition}

For the proof of Proposition \ref{estim:simples} we will use the following lemma:
\begin{lemma}\label{estim:simples:fraca}
Let $x_i\to x_0$ be an isolated simple blow-up point for $\{u_i\}$ and let $\rho$ be small. Then there exist $C,\delta>0$ such that
\begin{equation}\notag
M_i^{\lambda_i}|\nabla^k u_i|(\psi_i(z))\leq C|z|^{2-k-n+\rho},
\end{equation}
for $z \in B^+_{\delta}(0)\backslash\{0\}$ and $k=0,1,2$. Here, $\lambda_i=(p_i-1)(n-2-\rho)-1$.
\end{lemma}

The proof of Lemma \ref{estim:simples:fraca} is analogous to Lemma 2.7 of \cite{ahmedou-felli1}. It uses the following maximum principle, which is Lemma A.2 of \cite{han-li}:
\begin{lemma}\label{pm:1}
Let $(N,g)$ be a Riemannian manifold and $\Omega\subset N$ be a connected open set with piecewise smooth boundary $\partial \Omega=\Gamma\cup\Sigma$. Let $h\in L^{\infty}(\Omega)$ and $\sigma\in L^{\infty}(\Sigma)$. Suppose that $u\in C^2(\Omega)\cap C^1(\bar{\Omega})$ , $u>0$ in $\bar{\Omega}$, satisfies
\begin{align}
\begin{cases}
\Delta_g u + hu\leq 0\,,&\text{in}\; \Omega\,,\notag
\\
\frac{\partial u}{\partial\nu} + \sigma u\leq 0\,,&\text{on}\; \Sigma\notag
\end{cases}
\end{align}
and $v\in C^2(\Omega)\cap C^1(\bar{\Omega})$ satisfies
\begin{align}
\begin{cases}
\Delta_g v + hv\leq 0\,,&\text{in}\; \Omega\,,\notag
\\
\frac{\partial v}{\partial\nu} + \sigma v\leq 0\,,&\text{on}\; \Sigma\,,\notag
\\
v\geq 0\,,&\text{on}\; \Gamma\,,\notag
\end{cases}
\end{align}
where $\nu$ denotes inward unit normal to $\Sigma$. Then $v\geq 0$ in $\bar{\Omega}$.
\end{lemma}
\begin{remark}\label{rk:estim:simples:fraca}
Suppose that $x_i\to x_0$ is an isolated simple blow-up point for $\{u_i\}$. Then, as a consequence of the estimates of Lemma \ref{estim:simples:fraca} and Proposition \ref{form:bolha}, we see that
there exists $C>0$ such that
$$
|\nabla^k v_i|(y)\leq CM_i^{\rho(p_i-1)}(1+|y|)^{2-k-n}
$$
for any $y\in B^+_{\delta M_i^{p_i-1}}(0)$ and $k=0,1,2$. 
\end{remark}


Now we are going to estimate $\tau_i=\frac{n}{n-2}-p_i$.
\begin{proposition}\label{lim:ei:tau}
Let $x_i\to x_0$ be an isolated simple blow-up point for $\{u_i\}$ and let $\rho>0$ be small. Then there exists $C>0$ such that
\ba\label{estim:lim:ei:tau}
\tau_i\leq
\begin{cases}
C\e_i^{1-2\rho+o_i(1)},&\text{for}\:n\geq 4,
\\
C\e_i^{1-2\rho+o_i(1)}\log(\e_i^{-1}) ,&\text{for}\:n=3.
\end{cases}
\end{align}
(Recall that we have set $\ei=M_i^{-(p_i-1)}$ in the beginning of Section \ref{sec:isolated:simple}.)
\end{proposition}
\bp
\bigskip
Let  $x_i\to x_0$ be an isolated simple blow-up point for the sequence $\{u_i\}$. 
In order to simplify our notations, we will omit the simbol $\psi_i$ in the rest of this proof. Hence, points $\psi_i(z)\in M$, for $z\in \Bei=\Bei(0)$, will be denoted simply by $z$.
In particular, $x_i=\psi_i(0)$ will be denoted by $0$ and $u_i\circ\psi_i$ by $u_i$.

We write the Pohozaev identity of Proposition \ref{Pohozaev} as
\begin{equation}\label{corol:pohoz:1}
P(u_i,r)=F_i(u_i,r)+\bar{F}_i(u_i,r)+\frac{\tau_i}{p_i+1}Q_i(u_i,r),
\end{equation}
for $r<\delta$, where
\\
$F_i(u,r)=-\int_{B_r^+}(z^b\partial_bu+\frac{n-2}{2}u)(L_{g_i}-\Delta)u\,dz$,
\\
$\bar{F}_i(u,r)=\frac{n-2}{2}\int_{\partial 'B_r^+}(\bar{z}^b\partial_bu+\frac{n-2}{2}u)h_{g_i}u\,d\bar{z}$,
\\
$Q_i(u,r)=
\frac{(n-2)^2}{2}\int_{\d 'B_r^+}f_i^{-\tau_i}u^{p_i+1}d\bar{z}
-(n-2)\int_{\d 'B_r^+}(\bar{z}^k\partial_kf_i)f_i^{-\tau_i-1}u^{p_i+1}d\bar{z}\,.
$

It follows from  Proposition \ref{form:bolha} that we can choose a subsequence such that
$$
\int_{\d 'B^+_{r_i}}u_i^{p_i+1}\geq c>0,
$$
where $r_i=R_i\ei\to 0$ and $R_i\to \infty$. Hence, for $r>0$ small, $Q_i(u_i,r)\geq c>0$.

Using the estimate of Lemma \ref{estim:simples:fraca} we obtain 
\begin{equation}\label{corol:pohoz:3}
P_i(u_i,r)\leq C\ei^{\frac{2\lambda_i}{p_i-1}}= C\ei^{n-2-2\rho+o_i(1)}.
\end{equation}
Changing variables,
$$
\bar{F}_i(u_i,r)
=\frac{n-2}{2}\ei^{-\frac{2}{p_i-1}+n-1}\int_{\d 'B^+_{r\ei^{-1}}}
\left(\bar{y}^b\partial_bv_i+\frac{n-2}{2}v_i\right)h_{g_i}(\epsilon_i \bar{y})v_i(\bar{y}) d\bar{y}.\notag
$$
Observe that $-\frac{2}{p_i-1}+n-2=-(n-2)\frac{\tau_i}{p_i-1}=o_i(1)$. By Remark \ref{rk:estim:simples:fraca},  
\begin{align}\label{corol:pohoz:4}
\bar{F}_i(u_i,r)&=\ei^{1-2\rho+o_i(1)}
\int_{\d 'B_{r\ei^{-1}}^+}
O((1+|\bar{y}|)^{2-n})O((1+|\bar{y}|)^{2-n})d\bar{y}\notag
\\
&\geq -C\ei^{1-2\rho+o_i(1)}\cdot
\begin{cases}
1,&\text{for}\:n\geq 4,
\\
\log(\ei^{-1}),&\text{for}\:n=3.
\end{cases}
\end{align}

Similarly,
\ba
F_i(u_i,r)
&=-\ei^{-\frac{2}{p_i-1}+n-2}
\int_{B_{r\epsilon_i^{-1}}^+}
(y^b\partial_bv_i+\frac{n-2}{2}v_i)
(L_{\hat{g}_i}-\Delta)v_idy\notag
\\
&=\ei^{-2\rho+o_i(1)}
\int_{B_{r\ei^{-1}}^+}
O((1+|y|)^{2-n})O(\ei|y|)O((1+|y|)^{-n})dy\notag
\\
&\geq -C\ei^{1-2\rho+o_i(1)}\cdot
\begin{cases}\notag
1,&\text{for}\:n\geq 4,
\\
\log(\ei^{-1}),&\text{for}\:n=3\,,
\end{cases}
\end{align}
where $(\hat{g}_i)_{kl}(y)=(g_i)_{kl}(\ei y)$.
This, together with the identities (\ref{corol:pohoz:1}), (\ref{corol:pohoz:3}), (\ref{corol:pohoz:4}) and the fact that $Q_i(u_i,r)\geq c>0$, gives the result.
\ep



Now, we are able to prove Proposition \ref{estim:simples}.
\begin{proof}[Proposition \ref{estim:simples}]
We will first need the following two claims.
\\\\
{\it{Claim 1}} Given a small $\sigma >0$, there exists $C>0$ such that
$$
\int_{\d 'B_{\sigma}^+}u_i(\psi_i(\bar{z}))^{p_i}d\bar{z}\leq CM_i^{-1} \,.
$$

\vspace{0.2cm}
If follows from Proposition \ref{form:bolha} that we can choose a subsequence such that 
$$
\int_{\d 'B_{r_i}^+}u_i^{p_i}(\psi_i(\bar{z}))d\bar{z}=M_i^{-(p_i-1)(n-1)+p_i}\int_{\d 'B_{R_i}^+}v_i(\bar{y})^{p_i}d\bar{y}\leq CM_i^{-1}\,.
$$
Here, $r_i=R_iM_i^{-(p_i-1)}\to 0$, $R_i\to\infty$ and we used Proposition \ref{lim:ei:tau} in the last inequality. On the other hand, by Lemma \ref{estim:simples:fraca},
\ba
\int_{\d 'B^+_{\sigma}\backslash\d 'B_{r_i}^+} u_i^{p_i}(\psi_i(\bar{z}))d\bar{z}
&\leq CM_i^{-\l_ip_i}\int_{\d 'B^+_{\sigma}\backslash\d 'B^+_{r_i}}|\bar{z}|^{(2-n+\rho)p_i}d\bar{z}\notag
\\
&\leq CM_i^{-\l_ip_i}(R_iM_i^{-(p_i-1)})^{(2-n-\rho)p_i+n-1}\leq o_i(1)M_i^{-1}\,.\notag
\end{align}
This proves Claim 1.
\\\\
{\it{Claim 2}} There exists $\sigma_1>0$ such that for all $0<\sigma<\sigma_1$ there exists $C=C(\sigma)$ such that  
$$
u_i(\psi_i(z))u_i(x_i)\leq C
$$ 
for any $z\in\d^+B^+_{\sigma}(0)$.

\vspace{0.2cm}
It is not difficult to see that if $\sigma_1>0$ is small we can find a conformal metric, still denoted by $g_i$, such that $R_{g_{i}}\equiv 0$ in $B^+_{\sigma_1}(0)$ and $h_{g_i}\equiv 0$ on $\d 'B^+_{\sigma_1}(0)$. 

We fix $\sigma\in (0, \sigma_1)$ and choose any $x_{\sigma}\in\psi_i(\d^+B^+_{\sigma}(0))$. 

If we set $w_i=u_i(x_{\sigma})^{-1}u_i\circ\psi_i$, then $w_i$ satisfies 
\ba\label{estim:simples:1}
\begin{cases}
\Delta_{g_i}w_i=0,&\text{in}\: B^+_{\sigma}(0)\,,
\\
\d_n w_i+(n-2)u_i(x_{\sigma})^{p_i-1}f_i^{-\tau_i}w_i^{p_i}=0,&\text{on}\: \d 'B^+_{\sigma}(0)\,.
\end{cases}
\end{align}
By the Harnack inequality of Lemma \ref{Harnack:han-li}, for each $\b>0$ there exists $C_{\b}>0$ such that
$$
C_{\b}^{-1}\leq w_i(z)\leq C_{\b}
$$
if $|z|>\b$. Observe that Lemma \ref{estim:simples:fraca} implies that $u_i(x_{\sigma})^{p_i-1}\to 0$ as $i\to\infty$. Hence, we can suppose that $w_i\to w$ 
in $C^2_{loc}(B^+_{\sigma}(0)\backslash\{0\})$ for some $w>0$ satisfying
\ba
\begin{cases}
\Delta_{g_0}w=0,&\text{in}\: B_{\sigma}^+(0)\backslash\{0\},
\\
\d_n w=0,&\text{on}\: \d 'B_{\sigma}^+(0)\backslash\{0\}.
\end{cases}
\end{align}
Here, $g_{0}=\lim_{i\to\infty} g_i$.
It follows from elliptic linear theory that 
$$
w(z)=aG(z)+b(z)\,\:\:\:\:\text{for}\: z\in B^+_{\sigma}(0)\backslash\{0\}\,,
$$
where $a\geq 0$. Here, $G$ is the Green's function so that
\ba
\begin{cases}\notag
\Delta_{g_0}G=0,&\text{in}\: B_{\sigma}^+(0)\backslash\{0\},
\\
G=0,&\text{on}\: \d^+B_{\sigma}^+(0),
\\
\d_n G=0,&\text{on}\: \d 'B_{\sigma}^+(0)\backslash\{0\},
\\
\lim_{|z|\to 0}|z|^{n-2}G(z)=1\,,
\end{cases}
\end{align}
and $b$ satisfies
\ba
\begin{cases}\notag
\Delta_{g_0}b=0,&\text{in}\: B_{\sigma}^+(0),
\\
\d_n b=0,&\text{on}\: \d 'B_{\sigma}^+(0)\,.
\end{cases}
\end{align}

We will prove that $a>0$. We set $r=|z|$.
Since the blow-up is isolated simple, $r\mapsto r^{\frac{1}{p_i-1}}\bar{u}_i(r)$ is decreasing in $(r_i,\sigma)$ (see Remark \ref{rk:def:equiv}). Taking the limit as $i\to \infty$, we conclude that $r\mapsto r^{\frac{n-2}{2}}\bar{w}(r)$ is decreasing in $(0,\sigma)$. Hence, $w$ has a non-removable singularity at the origin. Therefore $a>0$.

Observe that there exists $c_1>0$ such that
\begin{equation}\label{estim:simples:3}
-\int_{\d^+B^+_{\sigma}}\frac{\d w}{\d r}d\sigma_{\sigma}>c_1.
\end{equation}

Integrating by parts the first equation of (\ref{estim:simples:1}) we obtain
\ba\label{estim:simples:4}
0=\int_{B^+_{\sigma}}\Delta_{g_0}w_idz
&=\int_{\d^+B^+_{\sigma}}\frac{\d w_i}{\d r}d\sigma_{\sigma}-\int_{\d 'B^+_{\sigma}}\d_n w_id\bar{z}\notag
\\
&=\int_{\d^+B^+_{\sigma}}\left(\frac{\d w}{\d r}+o_i(1)\right)d\sigma_{\sigma}
+(n-2)u_i(x_{\sigma})^{-1}\int_{\d 'B^+_{\sigma}}f_i^{-\tau_i}(u_i\circ\psi_i)^{p_i}d\bar{z}\notag
\\
&\leq -\frac{c_1}{2} +Cu_i(x_{\sigma})^{-1}u_i(x_i)^{-1},
\end{align}
where we used the estimate (\ref{estim:simples:3}) and Claim 1 in the last inequality. 
This proves Claim 2.

Now we are going to prove the item (a). Suppose by contradiction it does not hold. Then passing to a subsequence we can choose $\{x'_i\}\subset M$ such that $d_{g_i}(x'_i,x_i)\to 0$ and 
\begin{equation}\label{estim:simples:2}
u_i(x_i)u_i(x'_i)|z'_i|^{n-2}\to\infty ,
\end{equation}
where $z'_i=\psi_i^{-1}(x'_i)$.

By Proposition \ref{form:bolha} we can assume that 
$R_iu_i(x_i)^{-(p_i-1)}\leq |z'_i|\leq\delta/2$ where $R_i\to\infty$. We set 
$v_i(y)=|z'_i|^{\frac{1}{p_i-1}}u_i(\psi_i(|z'_i|\,y))$
for $y\in B^+_{\delta |z'_i|^{-1}}(0)$. Hence, the origin is an isolated simple blow-up point for $\{v_i\}$. Thus, by Claim 2, there exists $C>0$ such that
$$
|z'_i|^{\frac{2}{p_i-1}}u_i(x_i)u_i(x'_i)=v_i(0)v_i(y'_i)\leq C
$$
where $y'_i=|z'_i|^{-1}z'_i$. This contradicts the hypothesis (\ref{estim:simples:2}).

Item (b) is just an application of Lemma \ref{pm:1}.
\ep
\begin{remark}\label{rk:estim:simples}
Suppose that $x_i\to x_0$ is an isolated simple blow-up point for $\{u_i\}$. Then, as a consequence of  Propositions \ref{form:bolha} and \ref{estim:simples}, we see that $v_i\leq CU$ in $B^+_{\delta M_i^{p_1-1}}(0)$.
\end{remark}

\section{The linearized equation}

In this section we will be interested in solutions of a certain type of linear problem. 
These solutions will be used in the blow-up estimates of the next section.

\bigskip\noindent
{\bf{Convention}} In this section, we will always use the conformal equivalence between $\Rn\cup\{\infty\}$ and $B^n$ realized by the inversion $F$ (see Section 2.2).

\bigskip
Let $r\mapsto 0\leq\chi(r)\leq 1$ be a smooth cut-off function such that $\chi(r)\equiv 1$ for $0\leq r\leq \delta$ and $\chi(r)\equiv 0$ for $r>2\delta$. 
We set $\chi_{\e}(r)=\chi(\e r)$.
Thus,  $\chi_{\e}(r)\equiv 1$ for $0\leq r\leq \delta\e^{-1}$ and $\chi_{\e}(r)\equiv 0$ for $r>2\delta\e^{-1}$.
\begin{proposition}\label{exist:phi}
Let $\{h_{kl}^{(i)}\}_{i=1}^{\infty}$, $k,l=1,...,n-1$, and $\{\ei\}_{i=1}^{\infty}$ be sequences. Suppose that $tr(h_{kl}^{(i)})=0$, for each $i$, and $0<\ei\to 0$, as $i\to\infty$. Then, for each $i$, there is a solution $\phi_i$ to
\ba
\begin{cases}\label{linear:8}
\Delta\phi_{i}(y)=-2\chi_{\ei}(|y|)\ei h_{kl}^{(i)}y_n(\d_k\d_l U)(y)\,,&\text{for}\:y\in\Rn\,,
\\
\d_n\phi_{i}(\bar{y})+nU^{\frac{2}{n-2}}\phi_{i}(\bar{y})=0\,,&\text{for}\:\bar{y}\in\d\Rn\,,
\end{cases}
\end{align} 
where $\Delta$ stands for the Euclidean Laplacian, satisfying 
\begin{equation}\label{estim:phi'}
|\nabla^r\phi_i|(y)\leq C\ei |h_{kl}^{(i)}|(1+|y|)^{3-r-n}\,,\:\:\:\:\text{for}\: y\in\Rn\,,\,r=0,1 \:\text{or}\:2\:, 
\end{equation} 
\begin{equation}\label{hip:phi}
\phi_i(0)=\frac{\d\phi_i}{\d y_1}(0)=...=\frac{\d\phi_{i}}{\d y_{n-1}}(0)=0
\end{equation}
and
\begin{equation}\label{orto-phi}
\int_{\d\Rn}U^{\frac{n}{n-2}}(\bar{y})\phi_i(\bar{y})\,d\bar{y}=0\,.
\end{equation}
\end{proposition}
\bp
We set
$$
f_i(F(y))=-2\chi_{\ei}(|y|)\ei h_{kl}^{(i)}y_n(\d_k\d_l U)(y)U^{-\frac{n+2}{n-2}}(y)\:\:\:\:\text{for}\:y\in\Rn\,. 
$$
Observe that $f_i$ can be extended as a smooth function to $B^n$. According to Lemma \ref{eigenvalues:Bn}, the coordinate functions $z_1,...,z_n$, taken with center $(0,...,0,-1/2)$, satisfy the equations (\ref{linear:3'}) with $\lambda=2$ and we also have
$z_j\circ F=\frac{-1}{n-2}U^{-1}J_j$, for $j=1,...,n-1$, and $z_n\circ F=\frac{1}{n-2}U^{-1}J_n$. Hence,
\begin{equation}\notag
\int_{B^n}f_i\,z_jdz=\frac{2}{n-2}\int_{\Rn}\chi_{\ei}(|y|)\ei h_{kl}^{(i)}y_n(\d_k\d_l U)(\d_jU)dy=0\,,\:\:\:\:j=1,...,n-1\,,
\end{equation}
and
$$
\int_{B^n}f_i\,z_ndz=\frac{-2}{n-2}\int_{\Rn}\chi_{\ei}(|y|)\ei h_{kl}^{(i)}y_n(\d_k\d_l U)\left(\frac{n-2}{2}U+y^b\d_b U\right)dy=0\,.
$$
Here, we used the fact that
\ba\label{rel:int:pol}
\int_{S_r^{n-2}}p_k=\frac{r^2}{k(k+n-3)}\int_{S_r^{n-2}}\Delta p_k
\end{align} 
for every homogeneous polynomial $p_k$ of degree $k$.
Thus, by elliptic linear theory, it is possible to find a smooth solution $\bar{\phi}_{\ei}$ to
\begin{equation}\label{linear:9}
\begin{cases}
\Delta \bar{\phi}_{\ei} = f_{i}\,,&\text{in}\:B^n\,,\\
\frac{\partial \bar{\phi}_{\ei}}{\partial \eta}+2\bar{\phi}_{\ei}=0\,,&\text{on}\:\partial B^n\,,
\end{cases}
\end{equation}
$L^2(B^n)$-ortogonal to the coordinate functions $z_1,...,z_n$. We recall that $\eta$ is the inward unit normal vector to $\d B^n$.

Let $G$ be the Green's function on $B^n$ so that 
\begin{equation}\notag
\Delta G(z,w)=\a_n\sum_{a=1}^{n}z_aw_a\,,\:\:\:\:\:\text{for}\: z\neq w\,,
\end{equation}
subject to the boundary condition $\left(\frac{\d}{\d\eta}+2\right)G=0$ on $\d B^n$. Here, $\a_n=|z_a|^{-2}_{L^2(B^n)}$, for any $a=1,...,n$, is a constant.
By the Green's formula, $G$ satisfies 
\ba
\varphi(z)=\sum_{a=1}^{n}\int_{B^n}\a_nz_aw_a\,\varphi(w)\,dw
&-\int_{B^n}G(z,w)\Delta\varphi(w)\,dw\notag
\\
&-\int_{\d B^n}G(z,w)\left(\frac{\d \varphi}{\d\eta}+2\varphi\right)(w)\,d\sigma(w)\notag
\end{align}
for any $\varphi\in C^2(B^n)$. In particular, $\bar{\phi}_{\ei}$ satisfies
$$
\bar{\phi}_{\ei}(z)=-\int_{B^n}G(z,w)f_i(w)dw\,.
$$
Therefore,
$$
|\bar{\phi}_{\ei}(z)|\leq C\ei |h_{kl}^{(i)}|\int_{B^n}|z-w|^{2-n}|w+(0,...,0,1)|^{-3}dw\,.
$$
It follows from the result in \cite{giraud}, p.150 (see also \cite{aubin}, p.108) that
$$\bar{\phi}_{\ei}(z)\leq C\ei |h_{kl}^{(i)}||z+(0,...,0,1)|^{-1}\leq C\ei |h_{kl}^{(i)}|(|F(z)|+1)\,.$$ 
Hence, $\phi_{\ei}=U(\bar{\phi}_{\ei}\circ F)$  satisfies the estimate (\ref{estim:phi'}).  
By the properties (\ref{propr:L}) and (\ref{propr:B}) of the operators $L_g$ and $B_g$, $\phi_{\ei}$ is a solution to the equations (\ref{linear:8}).

Now, we choose coefficients $c_{j,i}=\frac{1}{n-2}\frac{\d\phi_{\ei}}{\d y_j}(0)$, $j=1,...,n-1$, and $c_{n,i}=-\frac{2}{n-2}\phi_{\ei}(0)$ and define
$$
\phi_i=\phi_{\ei}+\sum_{a=1}^{n}c_{a,i}J_a\,.
$$
Then $\phi_i$ is also a solution to the equations (\ref{linear:8}) and satisfies the identity (\ref{hip:phi}). Since $\phi_{\ei}$ satisfies the estimate (\ref{estim:phi'}), we see that $|c_{a,i}|\leq C|h_{kl}^{(i)}|\ei$ for $a=1,...,n$. Hence, $\phi_i$ also satisfies the estimate (\ref{estim:phi'}).

Let us prove the identity (\ref{orto-phi}). Observe that $\bar{\phi}_i=(U^{-1}\phi_i)\circ F^{-1}$ also satisfies the equations (\ref{linear:9}) and $f_i$ is $L^2(B^n)$-ortogonal to the constant function $1$. 
Hence, integrating by parts the equations (\ref{linear:9}) we see that $\bar{\phi}_i$ is $L^2(\d B^n)$-ortogonal to the funcion $1$. This is the identity (\ref{orto-phi}).
\ep


The following result is an important estimate that will be used in the subsequent local blow-up analysis.
\begin{proposition}\label{posit:compl}
Let $\phi_i$, $h_{kl}^{(i)}$ and $\ei$ be as in Proposition \ref{exist:phi} and suppose that $n\geq 5$. Then $\phi_i$ satisfies
\begin{align}
-&\int_{\Bei(0)}\Jzi \ei h_{kl}^{(i)}y_n\d_k\d_lU\,dy\notag
\\
&-\int_{\Bei(0)}\JU \ei h_{kl}^{(i)}y_n\d_k\d_l\phi_i\,dy
\geq -C(n)|h_{kl}^{(i)}|^2\ei^{n-2}\delta^{2-n}\,.\notag
\end{align}
\end{proposition}
\begin{proof}
We first recall that we have denoted $\Bei=\Bei(0)\subset\Rn$ for short (see Section \ref{subsec:notations}).
Integrating by parts,
\ba\label{posit_compl1}
-&\int_{\Bei}\Jzi(\ei h_{kl}^{(i)}y_n\d_k\d_lU)dy\notag
\\
&\geq\int_{\Bei}\ei h_{kl}^{(i)}y_n\d_k\phi_i\d_lUdy
+\int_{\Bei}\ei h_{kl}^{(i)}y_ny_b\d_b\d_k\phi_i\d_lUdy\notag
\\
&\hspace{1cm}+\frac{n-2}{2}\int_{\Bei}\ei h_{kl}^{(i)}y_n\d_k\phi_i\d_lUdy
-C|h_{kl}^{(i)}|^2\ei^{n-2}\delta^{2-n}
\end{align}
and
\ba\label{posit_compl2}
-&\int_{\Bei}\JU(\ei h_{kl}^{(i)}y_n\d_k\d_l\phi_i)dy\notag
\\
&\geq\int_{\Bei}\ei h_{kl}^{(i)}y_n\d_kU\d_l\phi_idy
+\int_{\Bei}\ei h_{kl}^{(i)}y_ny_b\d_b\d_kU\d_l\phi_idy\notag
\\
&\hspace{1cm}+\frac{n-2}{2}\int_{\Bei}\ei h_{kl}^{(i)}y_n\d_kU\d_l\phi_idy
-C|h_{kl}^{(i)}|^2\ei^{n-2}\delta^{2-n}\,.
\end{align}
Here, the terms $C|h_{kl}^{(i)}|^2\ei^{n-2}\delta^{2-n}$ come from estimating the integrals over $\d^+\Bei$ using the estimate (\ref{estim:phi'}). Another integration by parts gives
\ba
\int_{\Bei}\ei h_{kl}^{(i)}y_ny_b&(\d_b\d_k\phi_i)\d_lUdy+\int_{\Bei}\ei h_{kl}^{(i)}y_ny_b(\d_b\d_kU)\d_l\phi_idy\notag
\\
&\geq-(n+1)\int_{\Bei}\ei h_{kl}^{(i)}y_n\d_k\phi_i\d_lUdy
-C|h_{kl}^{(i)}|^2\ei^{n-2}\delta^{2-n}\,.\notag
\end{align}
This, together with the inequalities \eqref{posit_compl1} and \eqref{posit_compl2}, gives
\begin{align}\notag
-&\int_{\Bei}\Jzi(\ei h_{kl}^{(i)}y_n\d_k\d_lU)dy\notag
\\
&\hspace{0.5cm}-\int_{\Bei}\JU(\ei h_{kl}^{(i)}y_n\d_k\d_l\phi_i)dy\notag
\\
&\hspace{1cm}\geq-\int_{\Bei}\ei h_{kl}^{(i)}y_n\d_k\phi_i\d_lUdy
-C|h_{kl}^{(i)}|^2\ei^{n-2}\delta^{2-n}\,.\notag
\end{align}
The result now follows from the following Claim:
\\\\
{\it{Claim}} $\:\:-\int_{\Bei}\ei h_{kl}^{(i)}y_n\d_k\phi_i\d_lUdy\geq -C|h_{kl}^{(i)}|^2\ei^{n-2}\delta^{2-n}$

\vspace{0.2cm}
Integrating by parts and using the first equation of (\ref{linear:8}),
\ba
-\int_{\Bei}\ei h_{kl}^{(i)}y_n\d_k\phi_i\d_lUdy
&\geq\int_{\Bei}\phi_i\ei h_{kl}^{(i)}y_n\d_k\d_lUdy
-C|h_{kl}^{(i)}|^2\ei^{n-2}\delta^{2-n}\notag
\\
&=-\frac{1}{2}\int_{\Bei}(\Delta \phi_i)\phi_idy
-C|h_{kl}^{(i)}|^2\ei^{n-2}\delta^{2-n}\,.\notag
\end{align}
It follows from the estimate (\ref{estim:phi'}) and the assumption over the dimension that
$$-\int_{\Bei}(\Delta \phi_i)\phi_idy\geq -\int_{\Rn}(\Delta \phi_i)\phi_idy
-C|h_{kl}^{(i)}|^2\ei^{n-2}\delta^{2-n}\,.$$ 
Hence, in order to prove the Claim, we will show that
\begin{equation}\label{posit_compl3}
-\int_{\Rn}(\Delta \phi_i)\phi_idy\geq 0
\end{equation}

We set $\bar{\phi}_i=(U^{-1}\phi_i)\circ F^{-1}$. Observe that, by the properties (\ref{propr:L}) and (\ref{propr:B}), $\bar{\phi}_i$ satisfies the equations (\ref{linear:9}) and we have
\begin{equation}\label{posit_compl4}
-\int_{\Rn}(\Delta \phi_i)\phi_idy
=-\int_{B^n}(\Delta_{B^n}\bar{\phi}_i)\bar{\phi}_idz\,.
\end{equation}
Now, integrating by parts in $B^n$, we obtain
\begin{equation}\label{posit_compl5}
-\int_{B^n}(\Delta_{B^n}\bar{\phi}_i)\bar{\phi}_idz
=\int_{B^n}|\nabla\bar{\phi}_i|^2_{B^n}dz
-2\int_{\d B^n}\bar{\phi}_i^2d\sigma\,,
\end{equation}

Using Lemma \ref{eigenvalues:Bn} we see that 
\[
\inf_{\bar{\phi}\in\mathcal{C}_1}\frac{\int_{B^n}|\nabla \bar{\phi}|^2dz}{\int_{\d B^n}\bar{\phi}^2d\sigma}=2\,,
\]
where $\mathcal{C}_1=\{\bar{\phi}\in H^1(B^n); \int_{\d B^n}\bar{\phi} d\sigma =0\}$. On the other hand, the identity  (\ref{orto-phi}) is equivalent to $\int_{\d B^n}\bar{\phi}_id\sigma=0$. Hence, 
\begin{equation}\label{posit_compl6}
\int_{B^n}|\nabla\bar{\phi}_i|^2_{B^n}dz
-2\int_{\d B^n}\bar{\phi}_i^2d\sigma\geq 0\,.
\end{equation}
Now the inequality (\ref{posit_compl3}) follows from the equalities \eqref{posit_compl4} and \eqref{posit_compl5} and the inequality \eqref{posit_compl6}. This proves the Claim.
\end{proof}

\section{Blow-up estimates}\label{sec:blowup:estim}

In this section we will give a pointwise estimate for a blow-up sequence $\{u_i\}$ in a neighborhood of an isolated simple blow-up point. The arguments given here are modifications of the ones given in \cite{khuri-marques-schoen} and \cite{marques} for the case of manifolds without boundary.

In what follows, we will make use of the notations $\ei=M_i^{-(p_i-1)}$, introduced in Section \ref{sec:isolated:simple}, $\tau_i=\frac{n}{n-2}-p_i$ and $u_i\in\mathcal{M}_i$, introduced in Section \ref{subsec:conformal_scalar}.

\bigskip\noindent
{\bf{Assumption}} In this section we assume that $n\geq 5$.

\bigskip
Let  $x_i\to x_0$ be an isolated simple blow-up point for the sequence $\{u_i\in\mathcal{M}_i\}$. We use conformal Fermi coordinates centered at $x_i$. Thus we will work with conformal metrics $\tilde{g}_i=\zeta_i^{\frac{4}{n-2}}g_i$ and sequences $\{\tilde{u}_i=\zeta_i^{-1}u_i\}$ and $\{\tilde{\e}_i\}$, where $\tilde{\e}_i=\tilde{u}_i(x_i)^{-(p_i-1)}=\ei$, since $\zeta_i(x_i)=1$. As observed in the Remark \ref{rk:mud_conf}, $x_i\to x_0$ is still an isolated blow-up point for the sequence $\{\tilde{u}_i\}$ and satisfies the same estimates of Proposition \ref{estim:simples} (since we have uniform control on the conformal factors $\zeta_i>0$, these estimates are preserved). Let $\fermilinha$ denote the $\tilde{g}_i$-Fermi coordinates centered at $x_i$. 

In order to simplify our notations, we will omit the simbols $\:\tilde{}\:$ and $\psi_i$ in the rest of this section. Thus, the metrics $\tilde{g}_i$ will be denoted by $g_i$ and points $\psi_i(z)\in M$, for $z\in B_{\delta'}^+(0)$, will be denoted simply by $z$.
In particular, $x_i=\psi_i(0)$ will be denoted by $0$ and $u_i\circ\psi_i$ by $u_i$. 

We set $v_i(y)=\ei^{\frac{1}{p_i-1}}u_i(\ei y)$ for $y\in \Beilinha=\Beilinha(0)$. We know that $v_i$ satisfies 
\begin{align}\label{eq:vi'}
\begin{cases}
L_{\hat{g}_i}v_i=0,&\text{in}\:\Beilinha,
\\
B_{\hat{g}_i}v_i+(n-2)\hat{f}_i^{-\tau_i}v_i^{p_i}=0,&\text{on}\:\d '\Beilinha,
\end{cases}
\end{align}
where $\hat{f}_i(y)=f_i(\ei y)$ and $\hat{g}_i$ is the metric with coefficients $(\hat{g}_i)_{kl}(y)=(g_i)_{kl}(\ei y)$.

Let $\phi_i$ be the solution to the linearized equation obtained in Proposition \ref{exist:phi} with $h_{kl}^{(i)}=h_{kl}(0)$, observing that the hypothesis $tr(h_{kl}^{(i)})=0$ is satisfied due to Proposition \ref{conf:fermi:thm}(iii).  
The main result of this section is the following:
\begin{proposition}\label{estim:blowup:compl}
There exist $C,\delta>0$ such that
\begin{align}
|v_i-(U+\phi_i)|(y)
&\leq C(|\d^2 g_i|+|\d g_i|^2)\epsilon_i^2(1+|y|)^{4-n}+C\ei^{n-3}(1+|y|)^{-1}\,,\notag
\\
|\nabla v_i-\nabla(U+\phi_i)|(y)
&\leq C(|\d^2 g_i|+|\d g_i|^2)\epsilon_i^2(1+|y|)^{3-n}+C\ei^{n-3}(1+|y|)^{-2}\,,\notag
\\
|\nabla^2 v_i-\nabla^2(U+\phi_i)|(y)
&\leq C(|\d^2 g_i|+|\d g_i|^2)\epsilon_i^2(1+|y|)^{2-n}+C\ei^{n-3}(1+|y|)^{-3}\,,\notag
\end{align}
for $|y|\leq\delta\ei^{-1}$. (See Section \ref{subsec:coodinate} for the notation $|\d^kg_i|$.)
\end{proposition}

In order to prove Proposition \ref{estim:blowup:compl} we will first prove some auxiliary results. 
\begin{lemma}\label{estim:blowup1}
There exist $\delta, C>0$ such that 

$$|v_i-U-\phi_i|(y)\leq C\max\{(|\d^2 g_i|+|\d g_i|^2)\ei^2, \ei^{n-3},\tau_i\}\,,$$
for $|y|\leq\delta\ei^{-1}$.
\end{lemma}
\begin{proof}
We consider $\delta<\delta'$ to be chosen later and set 
$$\Lambda_i=\max_{|y|\leq \delta\ei^{-1}} |v_i-U-\phi_i|(y)=|v_i-U-\phi_i|(y_i)\,,$$ 
for some $|y_i|\leq \delta\ei^{-1}$. From Remark \ref{rk:estim:simples}
we know that $v_i(y)\leq CU(y)$ for $|y|\leq \delta\ei^{-1}$. Hence, if there exists $c>0$ such that $|y_i|\geq c\ei^{-1}$, then
$$\Lambda_i=|v_i-U-\phi_i|(y_i)\leq C\,|y_i|^{2-n}\leq C\,\ei^{n-2}$$
where we used the estimate \eqref{estim:phi'} in the first inequality. This 
implies the stronger inequality $|v_i-U-\phi_i|(y)\leq C\,\epsilon_i^{n-2}$, for $|y|\leq \delta\ei^{-1}$. Hence, we can suppose that $|y_i|\leq \delta\ei^{-1}/2$. 

Suppose, by contradiction, the result is false. 
Then, choosing a subsequence if necessary, we can suppose that
\begin{equation}
\label{hipLambda}
\Lambda_i^{-1} (|\d^2 g_i|+|\d g_i|^2)\ei^2\to 0,\:\:\Lambda_i^{-1}\ei^{n-3}\:\text{and}\:\: \Lambda_i^{-1} \tau_i\to 0\,.
\end{equation}
We define
$$w_i(y)=\Lambda_i^{-1}(v_i-U-\phi_i)(y)\,,\:\:\:\:\text{for}\:\: |y|\leq \delta\ei^{-1}\,.$$
By the equations \eqref{eq:U} and \eqref{eq:vi'}, $w_i$ satisfies
\begin{equation}\label{wi}
\begin{cases}
L_{\hat{g}_i}w_i=Q_i\,,&\text{in}\:\Bei\,,\\
B_{\hat{g}_i}w_i+b_iw_i=\overline{Q}_i\,,&\text{on}\:\partial '\Bei\,,
\end{cases}
\end{equation}
where 

\begin{flushleft}
$b_i=(n-2)\hat{f}_i^{-\tau_i}\frac{v_i^{p_i}-(U+\phi_i)^{p_i}}{v_i-(U+\phi_i)}$,\\
$Q_i=-\Lambda_i^{-1}\left\{(L_{\hat{g}_i}-\Delta)(U+\phi_i)+\Delta \phi_i\right\}$,\\
$\overline{Q}_i=-\Lambda_i^{-1}\left\{(n-2)\hat{f}_i^{-\tau_i}(U+\phi_i)^{p_i}
-(n-2)U^{\frac{n}{n-2}}-nU^{\frac{2}{n-2}}\phi_i-\frac{n-2}{2}h_{\hat{g}_i}(U+\phi_i)\right\}$.
\end{flushleft}

Observe that, for any funcion $u$,
\begin{align}
(L_{\hat{g}_i}-\Delta)u(y)&=(\hat{g}^{kl}_i-\delta^{kl})(y)\d_k\d_lu(y)
+(\d_k\hat{g}^{kl}_i)(y)\d_lu(y)\notag
\\
&\hspace{2cm}-\frac{n-2}{4(n-1)}R_{\hat{g}_i}(y)u(y)
+\frac{\d_k \sqrt{\det \hat{g}_i}}{\sqrt{\det \hat{g}_i}}\hat{g}_i^{kl}(y)\d_lu(y)\notag
\\
&=(g^{kl}_i-\delta^{kl})(\epsilon_i y)\d_k\d_lu(y)
+\epsilon_i (\d_kg^{kl}_i)(\epsilon_i y)\d_lu(y)\notag
\\
&\hspace{2cm}-\frac{n-2}{4(n-1)}\epsilon_i^2R_{g_i}(\epsilon_i y)u(y)
+O(\ei^N|y|^{N-1})\d_lu(y)\,,\notag
\end{align}
where $N$ can be taken arbitrarily large since we are using conformal Fermi coordinates.
Hence, setting $N=n-3$,
\begin{align}\label{Qi}
Q_i(y)
&=-\Lambda_i^{-1}(g^{kl}_i-\delta^{kl})(\epsilon_i y)\d_{k}\d_l(U+\phi_i)(y)
-\Lambda_i^{-1}\epsilon_i (\d_kg^{kl}_i)(\epsilon_i y)\d_l(U+\phi_i)(y)\notag
\\
&\hspace{0.5cm}+\frac{n-2}{4(n-1)}\Lambda_i^{-1}\epsilon_i^2R_{g_i}(\epsilon_i y)(U+\phi_i)(y)-\Lambda_i^{-1}\Delta \phi_i(y)
+O(\Lambda_i^{-1}\ei^{n-3}|y|^{n-4}(1+|y|)^{1-n})\notag
\\
&=O\left(\Lambda_i^{-1}(|\d^2 g_i|+|\d g_i|^2)\ei^2(1+|y|)^{2-n}\right)
+O(\Lambda_i^{-1}\ei^{n-3}(1+|y|)^{-3})\,,
\end{align}
where we have used the identities \eqref{linear:8} and \eqref{estim:phi'} and Proposition \ref{exp:g}.

Observe that
\ba
&(n-2)\hat{f}_i^{-\tau_i}(U+\phi_i)^{p_i}
-(n-2)U^{\frac{n}{n-2}}-nU^{\frac{2}{n-2}}\phi_i\notag
\\
&\hspace{1cm}=(n-2)\left(\hat{f}_i^{-\tau_i}(U+\phi_i)^{p_i}-(U+\phi_i)^{\frac{n}{n-2}}\right)
+O(U^{\frac{4-n}{n-2}}\phi_i^2)\notag
\\
&\hspace{1cm}=(n-2)\hat{f}_i^{-\tau_i}\left((U+\phi_i)^{p_i}-(U+\phi_i)^{\frac{n}{n-2}}\right)\notag
\\
&\hspace{2cm}+(n-2)(\hat{f}_i^{-\tau_i}-1)(U+\phi_i)^{\frac{n}{n-2}}+O(U^{\frac{4-n}{n-2}}\phi_i^2)\,.\notag
\end{align}
Using

\begin{flushleft}
$U^{\frac{4-n}{n-2}}\phi_i^2=O(\ei^2|h_{kl}(0)|^2(1+|y|)^{2-n})$,\\
$h_{\hat{g}_i}(U+\phi_i)=O(\ei^2|\d^2 g_i|(1+|y|)^{3-n})$,\\
$\hat{f}_i^{-\tau_i}\left((U+\phi_i)^{p_i}-(U+\phi_i)^{\frac{n}{n-2}}\right)
=O(\tau_i(U+\phi_i)^{\frac{n}{n-2}}\log(U+\phi_i))=O(\tau_i(1+|y|)^{1-n})$,\\
$(\hat{f}_i^{-\tau_i}-1)(U+\phi_i)^{\frac{n}{n-2}}=O(\tau_i\log(f_i)(U+\phi_i)^{\frac{n}{n-2}})
=O(\tau_i(1+|y|)^{-n})$,
\end{flushleft}
where in the second line we used Proposition \ref{conf:fermi:thm}(iii),
we obtain
\begin{equation}\label{barQi}
\bar{Q}_i(\bar{y})= O\left(\Lambda_i^{-1}\ei^2(|\d^2 g_i|+|\d g_i|^2)(1+|\bar{y}|)^{3-n}\right)+
O\left(\Lambda_i^{-1}\tau_i(1+|\bar{y}|)^{1-n}\right)\,.
\end{equation}

Moreover,
\begin{equation}
\label{lim:bi}
b_i(y)\to nU^{\frac{2}{n-2}}\,,\:\:\:\:\text{in}\: C^2_{loc}(\Rn)\,,
\end{equation}
and
\begin{equation}
\label{estim:bi}
b_i(y)\leq C(1+|y|)^{-2}\,,\:\:\:\:\text{for}\: |y|\leq\delta\ei^{-1}\,. 
\end{equation}

Since $|w_i|\leq |w_i(y_i)|=1$, we can use standard elliptic estimates to conclude that $w_i\to w$, in $C_{loc}^2(\mathbb{R}_+^n)$, for some function $w$, choosing a subsequence if necessary. From the identities \eqref{hipLambda}, \eqref{Qi}, \eqref{barQi} and  \eqref{lim:bi}, we see that $w$ satisfies
\begin{equation}\label{w}
\begin{cases}
\Delta w=0\,,&\text{in}\:\mathbb{R}_+^n\,,\\
\d_n w+nU^{\frac{2}{n-2}}w=0\,,&\text{on}\:\partial\mathbb{R}_+^n\,.
\end{cases}
\end{equation}

\bigskip
\noindent
{\it{Claim}} $\:\:w(y)=O((1+|y|)^{-1})$, for $y\in \Rn$.

\vspace{0.2cm}
Choosing $\delta>0$ sufficiently small, we can consider the Green's function $G_i$ for the conformal Laplacian $L_{\hat{g}_i}$ in $\Bei$ subject to the boundary conditions $B_{\hat{g}_i} G_i=0$ on $\partial'\Bei$ and $G_i=0$ on $\partial^+\Bei$. Let $\eta_i$ be the inward unit normal vector to $\partial^+\Bei$. Then the Green's formula gives
\begin{align}\label{wG}
w_i(y)&=-\int_{\Bei}G_i(\xi,y)Q_i(\xi) \,dv_{\hat{g}_i}(\xi)
+\int_{\partial^+\Bei}\frac{\partial G_i}{\partial\eta_i}(\xi,y)w_i(\xi)\,d\sigma_{\hat{g}_i}(\xi)\notag
\\
&\hspace{1cm}
+\int_{\partial'\Bei}G_i(\xi,y)\left(b_i(\xi)w_i(\xi)-\overline{Q}_i(\xi)\right)\,d\sigma_{\hat{g}_i}(\xi)\,.
\end{align}
Using the estimates \eqref{Qi}, \eqref{barQi} and \eqref{estim:bi}  in the equation \eqref{wG}, we obtain
\begin{align}
|w_i(y)|
\leq &\:C\Lambda_i^{-1}(|\d^2 g_i|+|\d g_i|^2)\ei^2\int_{\Bei}|\xi-y|^{2-n}(1+|\xi|)^{2-n}d\xi\notag
\\
&\hspace{0.5cm}+C\Lambda_i^{-1}\ei^{n-3}\int_{\Bei}|\xi-y|^{2-n}(1+|\xi|)^{-3}d\xi
+C\int_{\d'\Bei}|\bar{\xi}-y|^{2-n}(1+|\bar{\xi}|)^{-2}d\bar{\xi}\notag
\\
&\hspace{1.5cm}+C\Lambda_i^{-1}(|\d^2 g_i|+|\d g_i|^2)\ei^2
\int_{\d'\Bei}|\bar{\xi}-y|^{2-n}(1+|\bar{\xi}|)^{3-n}d\bar{\xi}\notag
\\
&\hspace{0.5cm}+C\Lambda_i^{-1}\tau_i\int_{\d'\Bei}|\bar{\xi}-y|^{2-n}(1+|\bar{\xi}|)^{1-n}d\bar{\xi}
+C\Lambda_i^{-1}\ei^{n-2}\int_{\d^+\Bei}|\xi-y|^{1-n}d\sigma(\xi)\,,\notag
\end{align}
for $|y|\leq \delta\ei^{-1}/2$. Here, we have used the fact that $|G_i(x,y)|\leq C\,|x-y|^{2-n}$ for $|y|\leq \delta\ei^{-1}/2$ and, since $v_i(y)\leq CU(y)$, $|w_i(y)|\leq C\Lambda_i^{-1}\ei^{n-2}$ for $|y|=\delta\ei^{-1}$. Hence, using Lemma \ref{estim:int} and the assumption over the dimension,
\begin{equation}\label{estim:wi}
|w_i(y)|\leq C\,\left((1+|y|)^{-1}+\Lambda_i^{-1}(|\d^2 g_i|+|\d g_i|^2)\ei^2
+\Lambda_i^{-1}\ei^{n-3}+\Lambda_i^{-1}\tau_i\right)
\end{equation}
for $|y|\leq \delta\ei^{-1}/2$.
The Claim now follows from the
hypothesis (\ref{hipLambda}).

Now, we can use the claim above and Lemma \ref{classifLinear} to see that 
$$w(y)=\sum_{j=1}^{n-1}c_j\d_jU(y)
+c_n\left(\frac{n-2}{2}U(y)+y^b\partial_b U(y)\right)\,,$$
for some constants $c_1,...,c_n$.
It follows from the identity (\ref{hip:phi}) that $w_i(0)=\frac{\partial w_i}{\partial y_j}(0)=0$ for $j=1,...,n-1$. Thus we conclude that $c_1=...=c_n=0$. Hence, $w\equiv 0$. Since $w_i(y_i)=1$, we have $|y_i|\to\infty$. This, together with the hypothesis (\ref{hipLambda}), contradicts the estimate (\ref{estim:wi}), since $|y_i|\leq \delta\ei^{-1}/2$, and concludes the proof of Lemma \ref{estim:blowup1}.
\end{proof}
\begin{lemma}\label{estim:tau}
There exists $C>0$ such that
$$\tau_i\leq C\max\{(|\d^2 g_i|+|\d g_i|^2)\ei^2,\ei^{n-3}\}\,.$$
\end{lemma}
\begin{proof}
Suppose, by contradiction, the result is false. Then we can suppose that
\begin{equation}\label{hip:tau:epsilon}
\tau_i^{-1}(|\d^2 g_i|+|\d g_i|^2)\ei^2,\:\:\tau_i^{-1}\ei^{n-3}\to 0
\end{equation}
and, by Lemma \ref{estim:blowup1}, there exists $C>0$ such that
$$|v_i-(U+\phi_i)|(y)\leq C\tau_i\,,\:\:\:\:\text{for}\:\:|y|\leq\delta\ei^{-1}\,.$$
We define
$$w_i(y)=\tau_i^{-1}(v_i-(U+\phi_i))(y)\,,\:\:\:\:\text{for}\:\: |y|\leq \delta\ei^{-1}\,. $$
Then $w_i$ satisfies the equations \eqref{wi} with 

\begin{flushleft}
$b_i=(n-2)\hat{f}_i^{-\tau_i}\frac{v_i^{p_i}-(U+\phi_i)^{p_i}}{v_i-(U+\phi_i)}$,\\
$Q_i=-\tau_i^{-1}\left\{(L_{\hat{g}_i}-\Delta)(U+\phi_i)+\Delta \phi_i\right\}$,\\
$\overline{Q}_i=-\tau_i^{-1}\left\{(n-2)\hat{f}_i^{-\tau_i}(U+\phi_i)^{p_i}
-(n-2)U^{\frac{n}{n-2}}-nU^{\frac{2}{n-2}}\phi_i-\frac{n-2}{2}h_{\hat{g}_i}(U+\phi_i)\right\}$.
\end{flushleft}
Similarly to the estimates \eqref{Qi} and \eqref{barQi}  we have
\begin{align}
|Q_i(y)|&\leq C\tau_i^{-1}(|\d^2 g_i|+|\d g_i|^2)\ei^2(1+|y|)^{2-n}
+C\tau_i^{-1}\ei^{n-3}(1+|y|)^{-3}\,,\label{Qi2}
\\
|\overline{Q}_i(y)|&\leq C\tau_i^{-1}(|\d^2 g_i|+|\d g_i|^2)\ei^2(1+|y|)^{3-n}+
C(1+|y|)^{1-n}\label{barQi2}
\end{align}
and $b_i$ satisfies the estimate \eqref{estim:bi}.

By definition, $w_i\leq C$ and, by elliptic standard estimates, we can suppose that $w_i\to w$ in $C^2_{loc}(\mathbb{R}^n_+)$ for some function $w$. By the identity (\ref{lim:bi}) and the estimates \eqref{Qi2} and \eqref{barQi2} we see that $w$ satisfies the equations \eqref{w}.
Recall that $J_n(y)=\frac{n-2}{2}U(y)+y^b\partial_b U(y)$ also satisfies the equations \eqref{w} (see Section 2.2). 

Let $\eta_i$ be the inward unit normal vector to $\partial^+\Bei$. Using the Green's formula, we have
\begin{align}\label{psiQi}
\int_{\partial'\Bei}J_n\cdot (B_{\hat{g}_i}w_i+b_i w_i)\,d\sigma_{\hat{g}_i}
&=\int_{\partial'\Bei}(B_{\hat{g}_i}(J_n)+b_iJ_n)\cdot w_i\,d\sigma_{\hat{g}_i}
\\
&\hspace{1cm}+\int_{\partial^+\Bei}\left(\frac{\partial J_n}{\partial \eta_i}w_i
-J_n\frac{\partial w_i}{\partial \eta_i}\right)\,d\sigma_{\hat{g}_i}\notag\\
&\hspace{1.5cm}+\int_{\Bei}\left(w_i L_{\hat{g}_i}(J_n)-J_n L_{\hat{g}_i}(w_i)\right) \,dv_{\hat{g}_i}\,.\notag
\end{align}
It follows from the estimate (\ref{estim:phi'}) and the hypothesis (\ref{hip:tau:epsilon}) that
\begin{equation}\label{psi:w}
\lim_{i\to\infty}\int_{\partial^+\Bei}\left(\frac{\partial J_n}{\partial\eta_i}w_i
-J_n\frac{\partial w_i}{\partial\eta_i}\right)\,d\sigma_{\hat{g}_i}= 0\,.
\end{equation}
Using the first equation of \eqref{wi}, the estimate \eqref{Qi2} and again the hypothesis (\ref{hip:tau:epsilon}), we have
\begin{equation}\label{psi:Q}
\lim_{i\to\infty}\int_{\Bei}J_n L_{\hat{g}_i}(w_i)\,dv_{\hat{g}_i}=\lim_{i\to\infty}\int_{\Bei}J_n Q_i\,dv_{\hat{g}_i} = 0\,.
\end{equation}

We will now derive a contradiction using the identity (\ref{psiQi}). 
First observe that
\begin{equation}\label{psi:r}
J_n(y)=\frac{n-2}{2}\frac{1-r^2}{(1+r^2)^{\frac{n}{2}}}\,,\:\:\:\:\text{if}\:\: y_n=0\,. 
\end{equation}
Here, $r^2=y_1^2+...+y_{n-1}^2$.
Then
\begin{align}
\int_{\d\Rn}J_n U^{\frac{n}{n-2}}d\bar{y}
&=\frac{n-2}{2}\sigma_{n-2}\int_{0}^{\infty}\frac{1-r^2}{(1+r^2)^n}r^{n-2}dr\notag
\\
&=\frac{n-2}{2}\sigma_{n-2}
\left(
\int_{0}^{1}\frac{1-r^2}{(1+r^2)^n}r^{n-2}dr
+\int_{1}^{\infty}\frac{1-r^2}{(1+r^2)^n}r^{n-2}dr
\right)=0\,,\notag
\end{align}
where in the last equality we change variables $s=r^{-1}$.
Now, observe that 
$$\lim_{i\to\infty}\tau_i^{-1}\left(\hat{f}_i^{-\tau_i}(y)(U+\phi_i)^{p_i}(y)-(U+\phi_i)^{\frac{n}{n-2}}(y) \right)=-\left(\log f(x_0)+\log U(y)\right)U^{\frac{n}{n-2}}(y)\,,$$
where $f=\lim_{i\to\infty}f_i$. Similarly to the estimate \eqref{barQi2}, we have
\ba
&\left|\overline{Q}_i(y)+
(n-2)\tau_i^{-1}\left(\hat{f}_i^{-\tau_i}(y)(U+\phi_i)^{p_i}(y)-(U+\phi_i)^{\frac{n}{n-2}}(y)\right)\right|\notag
\\
&\hspace{3cm}\leq C\tau_i^{-1}(|\d^2 g_i|+|\d g_i|^2)\epsilon_i^{2}(1+|y|)^{3-n}\,.\notag
\end{align}
Therefore, since $\int_{\partial\mathbb{R}^n_+}J_n U^{\frac{n}{n-2}}\,d\bar{y} =0$,
\begin{equation}\label{psi:barQ}
\lim_{i\to\infty}\int_{\partial'\Bei} J_n \bar{Q}_i \,d\sigma_{\hat{g}_i}
=(n-2)\int_{\partial\mathbb{R}_+^n}J_n\log (U)U^{\frac{n}{n-2}}\,d\bar{y}\,,
\end{equation}
where we have used the hypothesis (\ref{hip:tau:epsilon}).
\\\\
{\it{Claim}} $\:\:\int_{\partial\mathbb{R}_+^n}J_n\log (U)U^{\frac{n}{n-2}}\,d\bar{y}>0$.

\vspace{0.2cm}
By the identity (\ref{psi:r}),
$$\int_{\partial\mathbb{R}_+^n}J_n(\log U)U^{\frac{n}{n-2}}\,d\bar{y}
=-\frac{(n-2)^2}{4}\sigma_{n-2}\int_{0}^{\infty}\frac{1-r^2}{(1+r^2)^n}\log (1+r^2)r^{n-2}dr\,.$$
Changing variables $s=r^{-1}$, we obtain
$$\int_{0}^{\infty}\frac{1-r^2}{(1+r^2)^n}\log (1+r^2)r^{n-2}dr=2\int_{1}^{\infty}\frac{1-r^2}{(1+r^2)^n}\log (r)r^{n-2}dr<0\,,$$
which concludes the proof of the Claim.

On the other hand, the equation \eqref{psiQi} together with the equations \eqref{wi}, \eqref{lim:bi}, \eqref{w}, \eqref{psi:w} and \eqref{psi:Q} gives
\begin{align}\label{psiQi2}
\lim_{i\to\infty}\int_{\partial'\Bei} J_n \bar{Q}_i \,d\sigma_{\hat{g}_i}
&=\lim_{i\to\infty}\int_{\partial'\Bei}w_i\cdot (B_{\hat{g}_i}(J_n)+b_iJ_n)\,d\sigma_{\hat{g}_i}
+\lim_{i\to\infty}\int_{\Bei}w_iL_{\hat{g}_i}(J_n)\,dv_{\hat{g}_i}\notag
\\
&=\int_{\partial\mathbb{R}^n_+} w\cdot \left(\frac{\partial J_n}{\partial y_n}
+nU^{\frac{2}{n-2}}J_n\right)\,d\bar{y}
+\int_{\mathbb{R}_+^n}w\Delta J_n\,dy=0\,.
\end{align}
Here, we have used the fact that, by the identity (\ref{psi:barQ}), this limit should be independent of $\delta>0$ arbitrarily small. By the previous claim, this contradicts the identity (\ref{psi:barQ}).  
\end{proof}
\begin{proposition}\label{estim:blowup}
There exist $C,\delta>0$ such that
\begin{equation}\notag
|v_i-(U+\phi_i)|(y)\leq C\max\{(|\d^2 g_i|+|\d g_i|^2)\ei^2, \ei^{n-3}\}\,,
\end{equation}
for $|y|\leq\delta\ei^{-1}$.
\end{proposition}
\bp
This result follows from Lemmas \ref{estim:blowup1} and \ref{estim:tau}.
\ep

Now, we are able to prove Proposition \ref{estim:blowup:compl}.
\begin{proof}[Proposition \ref{estim:blowup:compl}]
We define
$$w_i(y)=(v_i-(U+\phi_i))(y)\,,\:\:\:\:\text{for}\:\: |y|\leq \delta\ei^{-1}\,. $$
Then $w_i$ satisfies the equations \eqref{wi} with

\begin{flushleft}
$b_i=(n-2)\hat{f}_i^{-\tau_i}\frac{v_i^{p_i}-(U+\phi_i)^{p_i}}{v_i-(U+\phi_i)}$,\\
$Q_i=-\left\{(L_{\hat{g}_i}-\Delta)(U+\phi_i)+\Delta \phi_i\right\}$,\\
$\overline{Q}_i=-\left\{(n-2)\hat{f}_i^{-\tau_i}(U+\phi_i)^{p_i}
-(n-2)U^{\frac{n}{n-2}}-nU^{\frac{2}{n-2}}\phi_i-\frac{n-2}{2}h_{\hat{g}_i}(U+\phi_i)\right\}$.
\end{flushleft}
Observe that $b_i$ satisfies the estimate \eqref{estim:bi}.
Similarly to the estimates \eqref{Qi}, \eqref{barQi} we have
\begin{align}
|Q_i(y)|&\leq C\ei^2(|\d^2 g_i|+|\d g_i|^2)(1+|y|)^{2-n}
+C\ei^{n-3}(1+|y|)^{-3}\,,\label{Qi3}
\\
|\overline{Q}_i(y)|&\leq C\ei^2(|\d^2 g_i|+|\d g_i|^2)(1+|y|)^{3-n}+
C\tau_i(1+|y|)^{1-n}\label{barQi3}\notag
\\
&\leq C\ei^2(|\d^2 g_i|+|\d g_i|^2)(1+|y|)^{3-n}
+C\ei^{n-3}(1+|y|)^{1-n}\,,
\end{align}
where in the last inequality we used Lemma \ref{estim:tau}.

The Green's formula gives
\begin{align}\label{wG3}
w_i(y)&=-\int_{\Bei}G_i(\xi,y)Q_i(\xi) \,dv_{\hat{g}_i}(\xi)
+\int_{\partial^+\Bei}\frac{\partial G_i}{\partial\eta_i}(\xi,y)w_i(\xi)\,d\sigma_{\hat{g}_i}(\xi)\notag
\\
&\hspace{1cm}+\int_{\partial'\Bei}G_i(\xi,y)\left(b_i(\xi)w_i(\xi)
-\overline{Q}_i(\xi)\right)\,d\sigma_{\hat{g}_i}(\xi)\,.
\end{align}
where $\eta_i$ is the inward unit normal vector to $\d^+\Bei$ and $G_i$ is the Green's function $G_i$ for the conformal Laplacian $L_{\hat{g}_i}$ in $\Bei$ subject to the boundary conditions $B_{\hat{g}_i} G_i=0$ on $\partial'\Bei$ and $G_i=0$ on $\partial^+\Bei$. Using the estimates \eqref{estim:bi}, \eqref{Qi3}, \eqref{barQi3} and Proposition \ref{estim:blowup} in equation \eqref{wG3}, as in the proof of Lemma \ref{estim:blowup1} we obtain
\begin{equation}\label{wi3}
|w_i(y)|\leq C\epsilon_i^{2}(|\d^2 g_i|+|\d g_i|^2)(1+|y|)^{-1}
+C\ei^{n-3}(1+|y|)^{-1}\,,
\end{equation}
for $|y|\leq\delta\ei^{-1}/2$. If $n=5$, we have the result. If $n\geq 6$, we plug the inequality \eqref{wi3} in the Green's formula \eqref{wG3} until we reach
\begin{equation}\notag
|w_i(y)|\leq C\ei^2(|\d^2 g_i|+|\d g_i|^2)(1+|y|)^{4-n}
+C\ei^{n-3}(1+|y|)^{-1}\,.
\end{equation}

The derivative estimates follow from elliptic theory, finishing the proof.
\ep


\section{Local blow-up analysis}

In this section we will prove the vanishing of the trace-free second fundamental form in an isolated simple blow-up point if $n\geq 7$. We will also prove a Pohozaev sign condition that will be used later in the study of the blow-up set. The basic tool here will be the Pohozaev-type the identity of Section 3 and the blow-up estimates of Section 6.

\subsection{Vanishing of the trace-free 2nd fundamental form}

The vanishing of $\pi_{kl}$, the trace-free 2nd fundamental form of the boundary, in an isolated simple blow-up point is stated as follows:
\begin{theorem}\label{anul:umb}
Suppose that $n\geq 7$. Let $x_i\to x_0$ be an isolated simple blow-up point for the sequence $\{u_i\in\mathcal{M}_i\}$. Then

$$
|\pi_{kl}(x_i)|^2\leq C\ei\,.
$$ 
In particular, 
$\pi_{kl}(x_0)=0$.
\end{theorem}
\bp  
Let  $x_i\to x_0$ be an isolated simple blow-up point for the sequence $\{u_i\}$. We use conformal Fermi coordinates centered at $x_i$. Thus we will work with conformal metrics $\tilde{g}_i=\zeta_i^{\frac{4}{n-2}}g_i$ and sequences $\{\tilde{u}_i=\zeta_i^{-1}u_i\}$ and $\{\tilde{\e}_i\}$, where $\tilde{\e}_i=\tilde{u}_i(x_i)^{-(p_i-1)}=\ei$, since $\zeta_i(x_i)=1$. As observed in the Remark \ref{rk:mud_conf}, $x_i\to x_0$ is still an isolated blow-up point for the sequence $\{\tilde{u}_i\}$ and satisfy the same estimates of Proposition \ref{estim:simples} (since we have uniform control on the conformal factors $\zeta_i>0$, these estimates are preserved). Let $\fermi$ denote the $\tilde{g}_i$-Fermi coordinates centered at $x_i$. 

In order to simplify our notations, we will omit the simbols $\:\tilde{}\:$ and $\psi_i$ in the rest of this section. Thus, the metrics $\tilde{g}_i$ will be denoted by $g_i$ and points $\psi_i(z)\in M$, for $z\in B_{\delta}^+(0)$, will be denoted simply by $z$.
In particular, $x_i=\psi_i(0)$ will be denoted by $0$ and $u_i\circ\psi_i$ by $u_i$. 

We set $v_i(y)=\ei^{\frac{1}{p_i-1}}u_i(\ei y)$ for $y\in \Bei=\Bei(0)$. We know that $v_i$ satisfies 
\begin{align}\notag
\begin{cases}
L_{\hat{g}_i}v_i=0,&\text{in}\:\Bei,
\\
B_{\hat{g}_i}v_i+(n-2)\hat{f}_i^{-\tau_i}v_i^{p_i}=0,&\text{on}\:\d '\Bei,
\end{cases}
\end{align}
where $\hat{f}_i(y)=f_i(\ei y)$ and $\hat{g}_i$ is the metric with coefficients $(\hat{g}_i)_{kl}(y)=(g_i)_{kl}(\ei y)$.
Observe that, from Remark \ref{rk:estim:simples}, we know that $v_i\leq CU$ in $B^+_{\delta\ei^{-1}}$.

By Proposition \ref{conf:fermi:thm} (iii), we can suppose that $h(0)=h_{,\,k}(0)=0$. In particular, $\pi_{kl}(0)=h_{kl}(0)$.
Recall that we use indices $1\leq k,l\leq n-1$ and $1\leq a,b \leq n$ when working with coordinates.
In many parts of the proof we will use the identity (\ref{rel:int:pol}).

We write the Pohozaev identity of Proposition \ref{Pohozaev} as
\begin{equation}\label{Pohoz}
P(u_i,r)=F_i(u_i,r)+\bar{F}_i(u_i,r)+\frac{\tau_i}{p_i+1}Q_i(u_i,r)\,,
\end{equation}
where
\\
$F_i(u,r)=-\int_{B_r^+}(z^b\partial_bu+\frac{n-2}{2}u)(L_{g_i}-\Delta)u\,dz$,
\\
$\bar{F}_i(u,r)=\frac{n-2}{2}\int_{\partial 'B_r^+}(\bar{z}^b\partial_bu+\frac{n-2}{2}u)h_{g_i}u\,d\bar{z}$,
\\
$Q_i(u,r)=
\frac{(n-2)^2}{2}\int_{\d 'B_r^+}f_i^{-\tau_i}u^{p_i+1}d\bar{z}
-(n-2)\int_{\d 'B_r^+}(\bar{z}^k\partial_kf)f_i^{-\tau_i-1}u^{p_i+1}d\bar{z}$.


We choose $r>0$ small enough such that 
$Q_i(u_i,r)\geq 0$. 
For the term $\bar{F}_i$ we have,
\begin{equation}
\bar{F}_i(u_i,r)=\frac{n-2}{2}\epsilon_i^{-\frac{2}{(p_i-1)}+n-2}\int_{\partial 'B_{r\epsilon_i^{-1}}^+}
\left(\bar{y}^b\partial_bv_i+\frac{n-2}{2}v_i\right)\ei h_{g_i}(\epsilon_i \bar{y})v_i(\bar{y}) d\bar{y}\,,\notag
\end{equation}

Since $h(0)=h_{,\,k}(0)=0$ and the fact that, according to Proposition \ref{lim:ei:tau}, 
$\lim_{i\to\infty}\epsilon_i^{-\frac{2}{p_i-1}+n-2}=\lim_{i\to\infty}\epsilon_i^{-(n-2)\frac{\tau_i}{p_i-1}}= 1$,
we have
\begin{align}\label{estim:barF}
\bar{F}_i(u_i,r)&=(1+o_i(1))\int_{\partial 'B_{r\epsilon_i^{-1}}^+}O((1+|\bar{y}|)^{2-n})O(\epsilon_i^{3}|\d^3g_i||\bar{y}|^2)O((1+|\bar{y}|)^{2-n})d\bar{y}\notag
\\
&\geq -C\epsilon_i^{3}|\d^3g_i|\int_{\partial 'B_{r\ei^{-1}}^+}(1+|\bar{y}|)^{6-2n}d\bar{y}\,.
\end{align}


We set $\check{U}_i(z)=\ei^{-\frac{1}{p_i-1}}(U+\phi_i)(\ei^{-1}z)$, where $\phi_i$ is as in Section \ref{sec:blowup:estim}. 
Using the facts that $g_i^{nn}\equiv 1$ and $g_i^{kn}\equiv 0$ in Fermi coordinates, we have
\begin{align}
F_i(u_i,r)=&-\int_{B_r^+}(z^b\partial_bu_i
+\frac{n-2}{2}u_i)(L_{g_i}-\Delta)u_idz\notag
\\
&=-\epsilon_i^{-\frac{2}{(p_i-1)}+n-2}\int_{B_{r\epsilon_i^{-1}}^+}(y^b\partial_bv_i+\frac{n-2}{2}v_i)\notag
(L_{\hat{g}_i}-\Delta)v_idy\,,
\notag
\end{align}

\begin{align}
F_i(\check{U}_i,r)
&=-\int_{B_r^+}(z^b\partial_b\check{U}_i
+\frac{n-2}{2}\check{U}_i)(L_{g_i}-\Delta)\check{U}_idz\notag
\\
&=-\epsilon_i^{-\frac{2}{(p_i-1)}+n-2}\int_{B_{r\epsilon_i^{-1}}^+}\left(y^b\partial_b(U+\phi_i)
+\frac{n-2}{2}(U+\phi_i)\right)\notag
(L_{\hat{g}_i}-\Delta)(U+\phi_i)dy\,.
\notag
\end{align}

It follows from Proposition \ref{estim:blowup:compl} that
\begin{align}\label{approx:F}
|F_i(u_i,r)-F_i(\check{U}_i,r)|
&\leq C\epsilon_i^{3}(|\d g_i|+|\d^2g_i|)(|\d^2g_i|+|\d g_i|^2)\int_{B_{r\epsilon_i^{-1}}^+}(1+|y|)^{5-2n}dy\notag
\\
&\hspace{1cm}+C\epsilon_i^{n-2}(|\d g_i|+|\d^2g_i|)\int_{B_{r\epsilon_i^{-1}}^+}(1+|y|)^{-n}dy\,.
\end{align}

\bigskip
We write
\begin{equation}\label{F:Ri}
F_i(\check{U}_i,r)=(1+o_i(1))\left\{R_i(U,U)+R_i(U,\phi_i)+R_i(\phi_i,U)+R_i(\phi_i,\phi_i)\right\}\,,
\end{equation}
where we have defined
$$R_i(w_1,w_2)=-\int_{B^+_{r\ei^{-1}}}(y^b\d_bw_1+\frac{n-2}{2}w_1)(L_{\hat{g}_i}-\Delta)w_2dy\,.$$

Using the identities (\ref{estim:barF}), (\ref{approx:F}) and  (\ref{F:Ri}) and the fact that $Q_i(u_i,r)\geq 0$ in the equality (\ref{Pohoz}), we have
\ba\label{P:Ri}
P(u_i,r)&\geq (1+o_i(1))\left\{R_i(U,U)+R_i(U,\phi_i)+R_i(\phi_i,U)+R_i(\phi_i,\phi_i)\right\}\notag
\\
&\hspace{0.5cm}-C(|\d g_i||\d^2g_i|+|\d g_i|^3+|\d^2g_i|^2+|\d^3g_i|)\,\ei^3\notag
\\
&\hspace{1cm}-C(|\d g_i|+|\d^2g_i|)\,\ei^{n-2}\log(\ei^{-1})\log r\,.
\end{align}
By Proposition \ref{exp:g} and the estimate (\ref{estim:phi'}),
\begin{align}\label{eq:Ri:U:phi}
R_i(U,\phi_i)+R_i(\phi_i,U)&=-\int_{B^+_{r\ei^{-1}}}\Jzi(L_{\hat{g}_i}-\Delta)Udy\notag
\\
&\hspace{0.5cm}-\int_{B^+_{r\ei^{-1}}}\JU(L_{\hat{g}_i}-\Delta)\phi_idy\notag
\\
&\geq-\int_{B^+_{r\ei^{-1}}}\Jzi(2\ei h_{kl}(0)y_n\d_k\d_lU)dy\notag
\\
&\hspace{0.5cm}-\int_{B^+_{r\ei^{-1}}}\JU(2\ei h_{kl}(0)y_n\d_k\d_l\phi_i)dy\notag
\\
&\hspace{0.5cm}-C\ei^3|h_{kl}(0)|(|\d^2g_i|+|\d g_i|^2)\int_{B^+_{r\ei^{-1}}}(1+|y|)^{5-2n}dy\,.\notag
\end{align}

Now we apply Proposition \ref{posit:compl} to this inequality to ensure that
\begin{equation}\label{estim:Ri:U:phi}
R_i(U,\phi_i)+R_i(\phi_i,U)\geq 
-C\left(\ei^3|h_{kl}(0)|(|\d^2g_i|+|\d g_i|^2)+|h_{kl}(0)|^2\ei^{n-2}r^{2-n}\right)\,.
\end{equation}

On the other hand, it follows from the estimate \eqref{estim:phi'} that 
\begin{equation}\label{estim:Ri:phi:phi}
R_i(\phi_i,\phi_i)=\ei^3|h_{kl}(0)|^2|\d g_i|\int_{B^+_{r\ei^{-1}}}O((1+|y|)^{5-2n})dy\,.
\end{equation}


We will now handle the term $R_i(U,U)$.
Observe that
\begin{align}
\partial_lU(y)&=-(n-2)\left((1+y_n)^2+|\bar{y}|^2 \right)^{-\frac{n}{2}}y_l\,,\notag
\\
\d_k\d_lU(y)
&=(n-2)\left((1+y_n)^2+|\bar{y}|^2 \right)^{-\frac{n+2}{2}}
\left(ny_k y_l-((1+y_n)^2+|\bar{y}|^2)\delta_{kl}\right)\,,\notag
\\
y^b\d_bU+\frac{n-2}{2} U
&=-\frac{n-2}{2}\left((1+y_n)^2+|\bar{y}|^2\right)^{-\frac{n}{2}}(|y|^2-1)\,.\notag
\end{align}
Using this we obtain
\begin{align}
R_i(U,U)&=\frac{(n-2)^2}{2}\int_{B^+_{r\epsilon_i^{-1}}}
\frac{|y|^2-1}{((1+y_n)^2+|y|^2)^{n+1}}\notag
\\
&\hspace{3cm}\cdot(g_i^{kl}-\delta^{kl})(\epsilon_i y)
\left(ny_ky_l-((1+y_n)^2+|\bar{y}|^2)\delta_{kl}\right)dy\notag
\\
&\hspace{0.5cm}-\frac{(n-2)^2}{2}\int_{B^+_{r\epsilon_i^{-1}}}
\frac{|y|^2-1}{((1+y_n)^2+|y|^2)^{n}}
\cdot \epsilon_i(\d_kg_i^{kl})(\epsilon_i y)y_ldy\notag
\\
&\hspace{0.5cm}-\frac{(n-2)^2}{8(n-1)}\int_{B^+_{r\epsilon_i^{-1}}}
\frac{|y|^2-1}{((1+y_n)^2+|y|^2)^{n-1}}
\cdot\epsilon_i^2 R_{g_i}(\epsilon_i y)dy\,.\notag
\end{align} 
Using Proposition \ref{exp:g}, we have
$$
R_i(U,U)\geq\frac{(n-2)^2}{2}(A_1+A_2+A_3+A_4)-C(|\d^2g_i|+|\d g_i|^2)\epsilon_i^{n-2}r^{2-n}\,,
$$
where

\begin{flushleft}
$\hspace{0.5cm}
A_1=n\int_{y_n=0}^{\infty}\int_{s=0}^{\infty}\frac{s^2+y_n^2-1}{(s^2+(y_n+1)^2)^{n+1}}
\left\{\int_{S_s^{n-2}}(g_i^{kl}-\delta^{kl})(\epsilon_i y)y_ky_l\,d\sigma_s(y)\right\}
dsdy_n$,
\\
$\hspace{0.5cm}
A_2=-\int_{y_n=0}^{\infty}\int_{s=0}^{\infty}\frac{s^2+y_n^2-1}{(s^2+(y_n+1)^2)^{n}}
\left\{\int_{S_s^{n-2}}(g_i^{kl}-\delta^{kl})(\epsilon_i y)\delta_{kl}\,d\sigma_s(y)\right\}
dsdy_n$,
\\
$\hspace{0.5cm}
A_3=-\int_{y_n=0}^{\infty}\int_{s=0}^{\infty}\frac{s^2+y_n^2-1}{(s^2+(y_n+1)^2)^{n}}
\left\{\epsilon_i\int_{S_s^{n-2}}(\d_kg_i^{kl})(\epsilon_i y)y_l\,d\sigma_s(y)\right\}
dsdy_n$,
\\
$\hspace{0.5cm}
A_4=\frac{-1}{4(n-1)}\int_{y_n=0}^{\infty}\int_{s=0}^{\infty}\frac{s^2+y_n^2-1}{(s^2+(y_n+1)^2)^{n-1}}
\left\{\epsilon_i^2\int_{S_s^{n-2}}R_{g_i}(\epsilon_i y)\,d\sigma_s(y)\right\}
dsdy_n$.
\end{flushleft}

Using Propositions \ref{exp:g} and \ref{conf:fermi:thm} we see that
\begin{align}
\int_{S_s^{n-2}}(g_i^{kl}-\delta_i^{kl})(\epsilon_i y)y_ky_l\,d\sigma_s
&=\sigma_{n-2}\ei^2\frac{y_n^2s^{n}}{n-1}\cdot 2|h_{kl}(0)|^2
+\ei^3|\d^3g_i|O(|(s,y_n)|^{n+3})\,,\notag
\\
\int_{S_s^{n-2}}(g_i^{kl}-\delta_i^{kl})(\epsilon_i y)\delta_{kl}\,d\sigma_s
&=\sigma_{n-2}\ei^2\cdot y_n^2s^{n-2}\cdot 2|h_{kl}(0)|^2\notag
\\
&\hspace{1cm}
+\ei^3|\d^3g_i|O(|(s,y_n)|^{n+1})\,,\notag
\\
\epsilon_i\cdot\int_{S_s^{n-2}}(\d_kg_i^{kl})(\epsilon_i y)y_l\,d\sigma_s
&=\ei^3|\d^3g_i|O(|(s,y_n)|^{n+1})\,,\notag
\\
\epsilon_i^2\cdot\int_{S_s^{n-2}}R_{g_i}(\epsilon_i y)\,d\sigma_s
&=-\sigma_{n-2}\ei^2\cdot s^{n-2}\cdot |h_{kl}(0)|^2\notag
\\
&\hspace{1cm}
+\ei^3(|\d^3g_i|+|\d^2g_i||\d g_i|)O(|(s,y_n)|^{n-1})\,,\notag
\end{align}
where in the last equality we used the fact that, by the Gauss equation, $R(0)+|h_{kl}(0)|^2=0$. We set $I=\int_{0}^{\infty}\frac{s^{n}}{(s^2+1)^{n}}ds\,$. 
Using Corollary \ref{int:s} and the four equalities above, we obtain
\begin{align}
A_1&=\sigma_{n-2}\epsilon_i^2\cdot
\frac{2n}{n-1}|h_{kl}(0)|^2
\int_{y_n=0}^{\infty}y_n^2
\left\{\int_{s=0}^{\infty}\frac{s^2+y_n^2-1}{(s^2+(y_n+1)^2)^{n+1}}s^{n}ds\right\}dy_n\notag
\\
&\hspace{2cm}+\epsilon_i^3|\d^3g_i|\int_{\Rn}O((1+|y|)^{5-2n})dy\notag
\\
&=\sigma_{n-2}\epsilon_i^2 I\cdot
\frac{n+1}{n-1} |h_{kl}(0)|^2
\int_{y_n=0}^{\infty}y_n^2(y_n+1)^{1-n}dy_n\notag
\\
&\hspace{0.5cm}+\sigma_{n-2}\epsilon_i^2 I\cdot |h_{kl}(0)|^2
\int_{y_n=0}^{\infty}y_n^2(y_n^2-1)(y_n+1)^{-1-n}dy_n\notag
\\
&\hspace{2cm}+\epsilon_i^3|\d^3g_i|\int_{\Rn}O((1+|y|)^{5-2n})dy\,,\notag
\end{align}
\begin{align}
A_2&=-\sigma_{n-2}\epsilon_i^2\cdot
2|h_{kl}(0)|^2
\int_{y_n=0}^{\infty}y_n^2
\left\{\int_{s=0}^{\infty}\frac{s^2+y_n^2-1}{(s^2+(y_n+1)^2)^{n}}s^{n-2}ds\right\}dy_n\notag
\\
&\hspace{2cm}+\epsilon_i^3|\d^3g_i|\int_{\Rn}O((1+|y|)^{5-2n})dy\notag
\\
&=-\sigma_{n-2}\epsilon_i^2 I\cdot 2|h_{kl}(0)|^2
\int_{y_n=0}^{\infty}y_n^2(y_n+1)^{1-n}dy_n\notag
\\
&\hspace{0.5cm}-\sigma_{n-2}\epsilon_i^2 I\cdot 2|h_{kl}(0)|^2
\int_{y_n=0}^{\infty}y_n^2(y_n^2-1)(y_n+1)^{-1-n}dy_n\notag
\\
&\hspace{2cm}+\epsilon_i^3|\d^3g_i|\int_{\Rn}O((1+|y|)^{5-2n})dy\,,\notag
\end{align}
\begin{equation}\notag
A_3=\epsilon_i^3|\d^3g_i|\int_{\Rn}O((1+|y|)^{5-2n})dy
\end{equation}
and
\begin{align}
A_4&=\sigma_{n-2}\epsilon_i^2\cdot
\frac{1}{4(n-1)}|h_{kl}(0)|^2
\int_{y_n=0}^{\infty}
\left\{\int_{s=0}^{\infty}\frac{s^2+y_n^2-1}{(s^2+(y_n+1)^2)^{n-1}}s^{n-2}ds\right\}dy_n\notag
\\
&\hspace{2cm}+\epsilon_i^3(|\d^3g_i|+|\d^2g_i||\d g_i|)\int_{\Rn}O((1+|y|)^{5-2n})dy\notag
\\
&=\sigma_{n-2}\epsilon_i^2 I\cdot 
\frac{1}{2(n-3)}|h_{kl}(0)|^2
\int_{y_n=0}^{\infty}(y_n+1)^{3-n}dy_n\notag
\\
&\hspace{0.5cm}+\sigma_{n-2}\epsilon_i^2 I\cdot 
\frac{1}{2(n-1)}|h_{kl}(0)|^2
\int_{y_n=0}^{\infty}(y_n^2-1)(y_n+1)^{1-n}dy_n\notag
\\
&\hspace{2cm}+\epsilon_i^3(|\d^3g_i|+|\d^2g_i||\d g_i|)\int_{\Rn}O((1+|y|)^{5-2n})dy\,.\notag
\end{align}

We set $I_k=\int_{0}^{\infty}\frac{y_n^k}{(1+y_n)^n}dy_n$. It follows from the above computations that
\ba\label{Ri:U:U1}
R_i&(U,U)\geq-C\ei^3(|\d^3g_i|+|\d^2g_i||\d g_i|)-C\ei^{n-2}r^{2-n}(|\d^2g_i|+|\d g_i|^2)\notag
\\
+&\sigma_{n-2}\ei^2I\cdot
\left\{\frac{n+1}{n-1}(I_3+I_2)+(I_3-I_2)-2(I_3+I_2)-2(I_3-I_2)\right\}|h_{kl}(0)|^2\notag
\\
+&\sigma_{n-2}\ei^2I\cdot
\left\{\frac{1}{2(n-3)}(I_3+3I_2+3I_1+I_0)+\frac{1}{2(n-1)}(I_3+I_2-I_1-I_0)\right\}|h_{kl}(0)|^2\notag
\\
&=\sigma_{n-2}\ei^2I\cdot (\a_3 I_3+\a_2 I_2+\a_1 I_1+\a_0 I_0)\cdot|h_{kl}(0)|^2\notag
\\
&\hspace{1cm}-C\ei^3(|\d^3g_i|+|\d^2g_i||\d g_i|)-C\ei^{n-2}r^{2-n}(|\d^2g_i|+|\d g_i|^2)\,,
\end{align}
where $\a_3=-2+\frac{1}{2(n-3)}+\frac{5}{2(n-1)}$, $\a_2=\frac{3}{2(n-3)}+\frac{5}{2(n-1)}$, 
$\a_1=\frac{3}{2(n-3)}-\frac{1}{2(n-1)}$ and $\a_0=\frac{1}{2(n-3)}-\frac{1}{2(n-1)}$.


By Lemma \ref{rel:int:t}, $I_2=\frac{n-4}{3}I_3$, $I_1=\frac{(n-4)(n-3)}{6}I_3$ and $I_0=\frac{(n-4)(n-3)(n-2)}{6}I_3$. Then a directy computation shows that
\begin{equation}\notag
\a_0 I_0+\a_1 I_1+\a_2 I_2+\a_3 I_3=\frac{n-6}{3}I_3\,.
\end{equation}
This, together with the inequality (\ref{Ri:U:U1}), implies that
\begin{align}\label{Ri:U:U2}
R_i(U,U)&\geq\sigma_{n-2}\ei^2\frac{n-6}{3}I\cdot I_3 |h_{kl}(0)|^2
-C\ei^{n-2}r^{2-n}(|\d^2g_i|+|\d g_i|^2)\notag
\\
&\hspace{1cm}-C\ei^3(|\d^3g_i|+|\d^2 g_i||\d g_i|)\,.
\end{align}

Hence, by the estimates (\ref{P:Ri}), (\ref{estim:Ri:U:phi}), (\ref{estim:Ri:phi:phi}) and (\ref{Ri:U:U2}),
\ba\label{P:J}
P(u_i,r)&\geq (1+o_i(1))\sigma_{n-2}\ei^2\frac{n-6}{3}I\cdot I_3|h_{kl}(0)|^2
-C\ei^{n-2}\log(\ei^{-1})r^{2-n}(|\d g_i|+|\d^2g_i|)\notag
\\
&\hspace{1cm}-C\ei^3(|\d^3g_i|+|\d g_i||\d^2g_i|+|\d^2g_i|^2+|\d g_i|^3)\,.
\end{align}
On the other hand, by Proposition \ref{estim:simples} we can assume that $\epsilon^{-\frac{1}{p_i-1}}_iu_i$ converges in $C^2_{loc}(B_{\delta}(0)\backslash\{0\})$ for $\delta>0$ small. Hence, for $r>0$ small fixed, $\ei^{-\frac{2}{p_i-1}}P(u_i,r)$ converges as $i\to\infty$ and 
\begin{equation}\label{estim:P}
P(u_i,r)\leq C\epsilon_i^{n-2}.
\end{equation}
Then the estimate (\ref{P:J}) together with the estimate (\ref{estim:P}) and our dimension assumption gives
$|h_{kl}(0)|^2\leq C\ei$.
This proves Theorem \ref{anul:umb}, since under our assumptions $\pi_{kl}(x_i)=h_{kl}(0)$.
\ep


\subsection{Pohozaev sign condition}

Now we will state and prove the Pohozaev sign condition.  

We set
$$P'(u,r)=\int_{\partial^+ B^+_r(0)}\left(\frac{n-2}{2}u\frac{\partial u}{\partial r}-\frac{r}{2}|\nabla u|^2+r\left|\frac{\partial u}{\partial r}\right|^2\right)d\sigma_r\,,$$
where $\nabla$ stands for the Euclidean gradient.
\begin{theorem}\label{cond:sinal}
Let $x_i\to x_0$ be a blow-up point for the sequence $\{u_i\in\mathcal{M}_i\}$. Assume that $\pi_{kl}(x_0)\neq 0$ and $n\geq 7$. 
We use Fermi coordinates $\fermi$ centered at $x_i$. 
For $0<\tau_i\to 0$, we set
$$
w_i(y)=\tau_i^{\frac{1}{p_i-1}}u_i(\psi_i(\tau_i y))\,,\:\:\:\:\:\text{for}\: y\in B_{\delta\tau_i^{-1}}^+(0)\,.
$$  
Suppose that the origin $0$ is an isolated simple blow-up point for the sequence $\{w_i\}$ and that $w_i(0)w_i\to G$ away from the origin, for some function $G$. Then
\begin{equation}\label{eq:cond:sinal}
\liminf_{r\to 0}P'(G,r)\geq 0\,.
\end{equation}
\end{theorem}
\bp
We will use conformal Fermi coordinates centered at $x_i$. Hence, we actually work with a sequence $\{\tilde{u}_i=\zeta_i^{-1}u_i\}$ and metrics $\tilde{g}_i=\zeta_i^{\frac{4}{n-2}}g_i$ and we have uniform control on the conformal factors $\zeta_i>0$. Since $\tau_i\to 0$ and $\zeta_i(x_i)=1$, we see that $\tilde{w}_i(0)\tilde{w}_i(y)\to G(y)$, where $\tilde{w}_i(y)=\tau_i^{\frac{1}{p_i-1}}\tilde{u}_i(\psi_i(\tau_i y))$. Thus, we will use the same notations and conventions of the proof of Theorem \ref{anul:umb}, omiting the symbols $\:\tilde{}\:$ and $\psi_i$.

Observe that  
$|\pi_{kl}(x_i)|\geq \frac12|\pi_{kl}(x_0)|$
for $i$ large. We will restringe our analysis to $B_{\check{\delta}}^+(0)\subset B_{\delta\tau_i^{-1}}^+(0)$, for some $\check{\delta}>0$ fixed. We set $\check{\e}_i=w_i(0)^{-(p_i-1)}\to 0$. Hence, $\check{\e}=\ei\tau_i^{-1}$. Let $\check{g}_i$ be the metric on $B_{\check{\delta}}^+(0)$ with coefficients $(\check{g}_i)_{kl}(y)=(g_i)_{kl}(\tau_i y)$ and denote by $\check{h}_{kl}$ the corresponding 2nd fundamental form.

Similarly to the estimate (\ref{P:J}), we have
\ba\label{P:J'}
P(w_i,r)&\geq (1+o_i(1))\sigma_{n-2}\check{\e}_i^2\frac{n-6}{3}I\cdot I_3|\check{h}_{kl}(0)|^2
-C\check{\e}_i^{n-2}\log(\check{\e}_i^{-1})r^{2-n}(|\d \check{g}_i|+|\d^2 \check{g}_i|)\notag
\\
&\hspace{1cm}-C\check{\e}_i^3(|\d^3\check{g}_i|+|\d \check{g}_i||\d^2\check{g}_i|+|\d^2\check{g}_i|^2+|\d \check{g}_i|^3)\,.
\end{align}

By the Young's inequality, 
$$
\check{\e}_i^{n-2}\log (\check{\e}_i^{-1})r^{2-n}|\d\check{g}_i|
\leq |\d\check{g}_i|^2\check{\e}_i^{n-2}\log(\check{\e}_i^{-1})^2r^{2-2n}+\check{\e}_i^{n-2}r^2\,.
$$
Hence, writing the inequality (\ref{P:J'}) in terms of the metric $g_i$ we have
\ba
P(w_i,r)&\geq (1+o_i(1))\sigma_{n-2}\e_i^2\frac{n-6}{3}I\cdot I_3|h_{kl}(0)|^2
-C\ei^2(|\d g_i|^2+|\d^2 g_i|)\check{\e}_i^{n-4}\log(\check{\e}_i^{-1})^2r^{2-2n}\notag
\\
&\hspace{1cm}-C\ei^3(|\d^3g_i|+|\d g_i||\d^2g_i|+\tau_i|\d^2g_i|^2+|\d g_i|^3)
-C\check{\e}_i^{n-2}r^{2}\notag
\\
&\geq-C\check{\e}_i^{n-2}r^{2}\,,\notag
\end{align}
for large $i$ and $r>0$ small fixed. Here, we used our dimension assumption and the fact that $|h_{kl}(0)|=|\pi_{kl}(x_i)|\geq \frac12|\pi_{kl}(x_0)|>0$ in the last inequality. Hence,  
$$
P'(G,r)=\lim_{i\to \infty}\check{\e}_i^{-\frac{2}{p_i-1}}P(w_i,r)\geq -Cr^2\,,
$$
where we also used Proposition \ref{lim:ei:tau}.
This proves Theorem \ref{cond:sinal}.
\ep



\section{Proof of Theorem \ref{compactness:thm'}}

In this section, we will prove Theorem \ref{compactness:thm'}.

The first proposition of this section states that every isolated blow-up point $x_i\to x_0$ is also simple, as long as $\pi_{kl}$, the boundary trace-free 2nd fundamental form, does not vanish at $x_0$. 
\begin{proposition}\label{isolado:impl:simples}
Let $x_i\to x_0$ be a blow-up point for the sequence $\{u_i\in\mathcal{M}_i\}$. Assume that $\pi_{kl}(x_0)\neq 0$ and $n\geq 7$. We use Fermi coordinates $\fermi$ centered at $x_i$. If $0<\tau_i\to 0$ or $\tau_i=1$, we set
$$
w_i(y)=\tau_i^{\frac{1}{p_i-1}}u_i(\psi_i(\tau_i y))\,,\:\:\:\:\:\text{for}\:y\in B_{\delta\tau_i^{-1}}^+(0)\,.
$$  
Suppose that the origin $0$ is an isolated blow-up point for the sequence  $\{w_i\}$. Then it is also isolated simple.
\end{proposition}
\bp
Suppose that the origin is an isolated blow-up point for $\{w_i\}$  but is not simple. By definition, passing to a subsequence, there are at least two critical points of $r\mapsto r^{\frac{1}{p_i-1}}\bar{w}_i(r)$ in an interval $(0,\bar{\rho}_i)$, $\bar{\rho}_i\to 0$. Let $r_i=R_iw_i(0)^{-(p_i-1)}\to ~0$ and $R_i\to \infty$ be as in Proposition \ref{form:bolha}. By Remark \ref{rk:def:equiv}, 
there is exactly one critical point in the interval $(0,r_i)$. Let $\rho_i$ be the second critical point. Then $\bar{\rho}_i>\rho_i\geq r_i$.

We set 
$ 
v_i(z)=\rho_i^{\frac{1}{p_i-1}}w_i(\rho_i z)
$,
for $z\in B^+_{\delta \rho_i^{-1}\tau_i^{-1}}(0)$. Observe that, since $\rho_i\geq r_i$,
$$
v_i(0)^{p_i-1}=\rho_i w_i(0)^{p_i-1}\geq R_i\to\infty\,.
$$
Hence, $v_i(0)\to\infty$. 

By the scaling invariance (see Remark \ref{scaling:inv}),
the origin is an isolated blow-up point for $\{v_i\}$. By the definitions, $r\mapsto r^{\frac{1}{p_i-1}}\bar{v}_i(r)$ has exactly one critical point in the interval $(0,1)$ and
\begin{equation}\label{isolado:impl:simples:1}
\frac{d}{dr}(r^{\frac{1}{p_i-1}}\bar{v}_i(r))|_{r=1}=0\,.
\end{equation}
Hence, the origin is an isolated simple blow-up point for $\{v_i\}$. It follows from Proposition \ref{estim:simples}(a) that $v_i(0)v_i$ is uniformly bounded in compact subsets of $\Rn\backslash\{0\}$. Using  the equations (\ref{eq:ui}), we can suppose that $v_i(0)v_i$ converges in $C^2_{loc}(\Rn\backslash\{0\})$ for some function $G$ satisfying 
\begin{equation*}
\begin{cases}
\Delta G=0\,,&\text{in}\:\Rn\backslash\{0\}\,,
\\
\d_nG=0\,,&\text{on}\:\d\Rn\backslash\{0\}\,.
\end{cases}
\end{equation*}
From elliptic linear theory we know that $G(z)=a|z|^{2-n}+b(z)$, where $b$ is harmonic on $\Rn$ with Neumann boundary condition on $\d\Rn$. It follows from Proposition \ref{estim:simples}(b) that $a>0$. Since $G>0$, $\liminf_{|z|\to\infty} b(z)\geq 0$. By the Liouville's theorem, $b$ is constant. By the equality (\ref{isolado:impl:simples:1}),
$$
\frac{d}{dr}(r^{\frac{n-2}{2}}h(r))|_{r=1}=0\,,
$$
which implies that $b=a>0$. This contradicts the sign condition of Theorem \ref{cond:sinal}.
\ep


The next proposition ensures that the set $\{x_1,...,x_N\}\subset\d M$ of points obtained in Proposition \ref{conj:isolados} can only contain isolated blow-up points for any blow-up sequence $\{u_i\}$ as long as $\pi_{kl}$ does not vanish at the blow-up point.
Recall that we denote by $D_{\delta}(x_0)$ the metric ball of $\d M$ with radius $\delta$ centered at $x_0\in\d M$ and by $\bar{g}$ the boundary metric.
\begin{proposition} \label{dist:unif}
Assume that $n\geq 7$.
Let $\b>0$ be small, $R>0$ be large and consider $C_0=C_0(\b,R)$ and $C_1=C_1(\b,R)$ as in Proposition \ref{conj:isolados}. Let $x_0\in\d M$ be a point such that $\pi_{kl}(x_0)\neq ~0$. 
Then there exists $\delta>0$ such that, for any $u\in\mathcal{M}_{p}$ satisfying $\max_{\d M}u\geq C_0$, the set $D_{\delta}(x_0)\cap\{x_1(u),...,x_N(u)\}$ consists of at most one point. Here, $x_1(u),...,x_N(u)\in\d M$, with $N=N(u)$, are the points obtained in Proposition \ref{conj:isolados}. 
\end{proposition}
\begin{proof}
Suppose the result is not true. Then there exist sequences $p_i\in\left(\frac{n}{n-2}-\beta,\frac{n}{n-2}\right]$ and $u_{i}\in\mathcal{M}_{p_i}$ with $\max_{\d M}u_i\geq C_0$, such that after relabeling the indices we have
$
x_{1}^{(i)},x_{2}^{(i)} \to x_0
$,
as $i \to \infty$. Here, we have set $x_{1}^{(i)}=x_1(u_{i}),...,x_{N_i}^{(i)}=x_{N_i}(u_{i})$ and $N_i=N(u_{i})$.

We define 
$$
s_i = d_{\bar{g}}(x_{1}^{(i)},x_{2}^{(i)})^{-\frac{1}{2}} \to \infty.
$$
\\\\
{\it{Claim 1}} There exist $1 \leq j_i \neq k_i \leq N_i$ such that
$
x_{j_i}^{(i)},x_{k_i}^{(i)} \in D_{2s_i^{-1}}(x_1^{(i)})
$,
$$
\sigma_i = d_{\bar{g}}(x_{j_i}^{(i)},x_{k_i}^{(i)}) \leq d_{\bar{g}}(x_{1}^{(i)},x_{2}^{(i)}),
$$ 
$$
d_{\bar{g}}(x_l^{(i)},x_m^{(i)}) \geq \frac12 \sigma_i\,,\:\:\:
\text{for all}\: x_l^{(i)},x_m^{(i)} \in D_{s_i\sigma_i}(x_{j_i}^{(i)}),\: l \neq m\,.
$$

\vspace{0.2cm}
Suppose that Claim 1 is false. Then there exist $x_{l_1}^{(i)},x_{m_1}^{(i)} \in D_{s_i^{-1}}(x_{1}^{(i)})$, $l_1 \neq m_1$,
with 
$$
\sigma_{1,i} = d_{\bar{g}}(x_{l_1}^{(i)},x_{m_1}^{(i)}) < \frac12 \sigma_{0,i}=\frac12 s_i^{-2}.
$$
If we repeat this procedure, we obtain sequences $x_{l_r}^{(i)},x_{m_r}^{(i)} \in D_{s_i\sigma_{r-1,i}}(x_{l_{r-1}}^{(i)})$, $l_r \neq m_r$, with
$$
\sigma_{r,i} = d_{\bar{g}}(x_{l_r}^{(i)},x_{m_r}^{(i)}) < \frac12\sigma_{r-1,i}.
$$
Since $N_i < \infty$, this procedure has to stop and we reach a contradiction. This proves Claim 1.

Using Claim 1 and a relabeling of indices, we find $x_{1}^{(i)},x_{2}^{(i)} \to x_0$ and $s_i \to \infty$ so that,
if $\sigma_i = d_{\bar{g}}(x_{1}^{(i)},x_{2}^{(i)})$, we have $s_i\sigma_i \to 0$ and 
$$
d_{\bar{g}}(x_l^{(i)},x_m^{(i)}) \geq \frac12 \sigma_i\,,\:\:\:
\text{for all}\: x_l^{(i)},x_m^{(i)} \in D_{s_i\sigma_i}(x_{1}^{(i)})\,,\: l \neq m\,.
$$
By the item (3) of Proposition \ref{conj:isolados}  we have
$u_{i}(x_{1}^{(i)}),u_{i}(x_{2}^{(i)})\rightarrow\infty$.

Now we use Fermi coordinates $\fermi$ centered at $x_{1}^{(i)}$ and set
$$
v_{i}(y)=\sigma_{i}^{\frac{1}{p_{i}-1}}u_{i}(\psi_i(\sigma_{i}y))\,,\:\:\:\:\:
\text{for}\: y\in B_{s_i}^+(0)\,.
$$

If
$x_{l}^{(i)}\in D_{s_i\sigma_i}(x_{1}^{(i)})$, we set $y_{l}^{(i)}=\sigma_{i}^{-1}\psi_i^{-1}(x_{l}^{(i)})\in \d 'B_{s_i}^+(0)$. In particular, $y_1^{(i)}=0$. Then each $y_{l}^{(i)}$ is a local
maximum of $v_{i}$ and by the item (3) of Proposition \ref{conj:isolados},
\begin{equation*}
\min_{l} \{|y-y_{l}^{(i)}|^{\frac{1}{p_{i}-1}}\}v_{i}(y)\leq C\,,\:\:\:\:\:
\text{for}\: y \in \d 'B_{\frac12 s_i}^+(0)\,.
\end{equation*}
Furthermore, $|y_{2}^{(i)}|=|y_{1}^{(i)}-y_{2}^{(i)}|=1$ and $\min_{l\neq m}|y_{l}^{(i)}-y_{m}^{(i)}|\geq \frac12+o_i(1)$. 
\\\\
{\it{Claim 2}} $v_{i}(y_{1}^{(i)}), v_{i}(y_{2}^{(i)}) \to \infty$.

\vspace{0.2cm}
If $v_{i}(y_{2}^{(i)})$ stays bounded but
$v_{i}(y_{1}^{(i)})\rightarrow\infty$, then $y_1^{(i)}=0$ is an isolated blow-up point for $\{v_i\}$ and hence is isolated simple. Since $v_{i}$ remains uniformly
bounded near $y_{2}^{(i)}$, it follows from Lemma \ref{Harnack:han-li} and Proposition \ref{estim:simples}  that $v_{i}(y_{2}^{(i)})\rightarrow 0$. This
is a contradiction since the item (1) of Proposition \ref{conj:isolados}  implies that
$$
\sigma_{i}\geq\max\{Ru_{i}(x_{1}^{(i)})^{-(p_{i}-1)},
Ru_{i}(x_{2}^{(i)})^{-(p_{i}-1)}\},
$$
thus 
\begin{equation}\label{dist:unif:1}
v_{i}(y_{1}^{(i)}),v_{i}(y_{2}^{(i)})\geq R^\frac{1}{p_i-1}.
\end{equation}
Of course the same argument holds if we  exchange the roles of $v_{i}(y_{1}^{(i)})$ and $v_{i}(y_{2}^{(i)})$. 

On the other hand, if both $v_{i}(y_{1}^{(i)})$ and $v_{i}(y_{2}^{(i)})$ remain
bounded, we can suppose that any other $v_{i}(y_{l}^{(i)})$ also does, using the same argument above. Then, after passing to a subsequence, $v_{i}\to v$ in $C^{2}_{loc}(\Rn)$ for some $v>0$ satisfying 
\begin{equation*}
\begin{cases}
\Delta v=0\,,&\text{in}\:\Rn\,,
\\
\d_nv+f(x_0)^{p_0-\frac{n}{n-2}}\,v^{p_0}=0\,,&\text{on}\:\d\Rn\,,
\end{cases}
\end{equation*}
and $\d_kv(0)=\d_kv(y_2)=0$ for $k=1,...,n-1$. Here $p_0=\lim_{i\to\infty}p_{i}\in [\frac{n}{n-2}-\b,\frac{n}{n-2}]$ and  $y_2 = \lim_{i \to \infty}y_{2}^{(i)}$.
Note that $|y_2|=1$. Then the Liouville-type theorems of \cite{hu} and \cite{li-zhu} yield that
$v\equiv 0$, which contradicts the inequalities (\ref{dist:unif:1}). This proves Claim 2.

It follows from Claim 2 that $0=y_{1}^{(i)}$ and $y_{2}^{(i)}$ are isolated
blow-up points for $\{v_i\}$. Thus Proposition \ref{isolado:impl:simples} implies that they are isolated simple.

Then, similarly to the proof of Proposition \ref{isolado:impl:simples},
\begin{equation*}
v_{i}(y_{1}^{(i)})v_{i}(y)\rightarrow
G(y)=a_{1}|y|^{2-n}+a_{2}|y-y_{2}|^{2-n}+b(y)
\end{equation*}
in $C^{2}_{loc}(\Rn\backslash S)$. Here, $S$ denotes the set of blow-up points for $\{v_{i}\}$,
$b(y)$ is a harmonic function on
$\Rn\backslash (S\backslash\{0,y_{2}\})$ with Neumann boundary condition and $a_{1}$,
$a_{2}>0$.  By the maximum principle, $b(y)\geq 0$. Hence, for $|y|$ near
0,
\begin{equation*}
G(y)=a_{1}|y|^{2-n}+b+O(|y|)
\end{equation*}
for some constant $b>0$.  This contradicts the sign condition of Theorem \ref{cond:sinal} and proves Proposition \ref{dist:unif}.
\end{proof}


Now we are able to prove Theorem \ref{compactness:thm'}.
\begin{proof}[Theorem \ref{compactness:thm'}]
Suppose by contradiction that $x_i\to x_0$ is a blow-up point for a sequence $\{u_i\in\mathcal{M}_{p_i}\}$ and $\pi_{kl}(x_0)\neq 0$. Let $x_{1}(u_i),...,x_{N(u_i)}(u_i)$ be the  points obtained in Proposition \ref{conj:isolados}. By the item (3) of this Proposition, we must have $d_g(x_i,x_{k_i}(u_i))\to 0$ for some $1\leq k_i\leq N(u_i)$. If $x_{k_i}=x_{k_i}(u_i)$,
it is  not difficult to see that $u_{i}(x_{k_i})\to \infty$. Thus $x_{k_i}\to x_0$ is a blow-up point for $\{u_i\}$. It follows from Propositions \ref{isolado:impl:simples} and \ref{dist:unif} that $x_{k_i}\to x_0$ is isolated simple. This contradicts Theorem \ref{anul:umb}.
\end{proof}


\section{Appendix}

In this section we will state some technical results that were used in the previous computations. 

Our first result is a modification of Proposition 2.7 in \cite{lee-parker}. The proof is similar.
\begin{lemma}\label{extensao:sol}
Let $(M,g)$ be a Riemannian manifold with boundary $\partial M$. Let $x\in \partial M$ and $\mathcal{U}\subset M$ be an open set containing $x$. 
Let $u$ be a weak solution to
\begin{equation}\notag
\begin{cases}
\Delta u = 0\,,&\text{in}\:\mathcal{U}\backslash\{x\}\,
\\
(\frac{\partial}{\partial \eta} +\psi)u=0\,,&\text{on}\:\mathcal{U}\cap\partial M\backslash\{x\}\,,
\end{cases}
\end{equation}
where $\eta$ is the inward unit normal vector to $\d M$. Suppose that $u\in L^{q}(\mathcal{U})$ for some $q>\frac{n}{n-2}$ and $u,\psi u\in L^1(\mathcal{U}\cap\partial M)$. 
Then $u$ is a weak solution to
\begin{equation}\notag
\begin{cases}
\Delta u = 0\,,&\text{in}\:\mathcal{U}\,,
\\
(\frac{\partial}{\partial \eta} +\psi)u=0\,,&\text{on}\:\mathcal{U}\cap\partial M\,.
\end{cases}
\end{equation}
\end{lemma}


The proof of the following lemma is similar to the result in \cite{giraud}, p.150 (see also \cite{aubin}, p.108).
\begin{lemma}\label{estim:int}
Let $\rho>0$ be small and suppose that $\rho\leq\beta\leq\beta+\rho\leq\alpha\leq n-\rho$. Then there exists $C=C(n,\rho)>0$ such that  
$$\int_{\R^n}|y-x|^{\beta -n}(1+|x|)^{-\alpha}dx\leq C(1+|y|)^{\beta-\alpha}$$
for any $y~\in~\R^{n+k}~\supset~\R^n$.
%
%
%
\end{lemma}

For the proof we decompose $\R^n$ in three regions
\\
$\mathcal{A}:=\{x\in\R^n;\,|x-y|\leq\frac{1}{2}|y|+\frac{1}{2} \}$,
\\
$\mathcal{B}:=\{x\in\R^n;\,|x-y|\geq\frac{1}{2}|y|+\frac{1}{2}, |x|\leq 2|y|+1 \}$,
\\
$\mathcal{C}:=\{x\in\R^n;\,|x|\geq 2|y|+1 \}$,
\\
and perform the estimates in each one separately.

\bigskip
The following Harnack-type inequality is Lemma A.1 of \cite{han-li}:
\begin{lemma}\label{Harnack:han-li}
Let $L$ be an operator of the form
$$
Lu=\d_a\left(\a^{ab}(x)\d_bu+\b^a(x)u\right)+\gamma^a(x)\d_au+\zeta(x)u\,,\:\:\:a,b=1,...,n
$$
and assume that for some constant $\Lambda>1$ the coefficient functions satisfy
$$
\Lambda^{-1}|\xi|^2\leq \a^{ab}(x)\xi_a\xi_b\leq\Lambda|\xi|^2,
$$
$$
|\b^a(x)|+|\gamma^a(x)|+|\zeta(x)|\leq \Lambda,
$$
for all $x\in B^+_3=B^+_3(0)$ and all $\xi\in \R^n$. If $|q(x)|\leq\Lambda$, for any $x\in \d' B^+_3$, and $u\in C^2(B^+_3\backslash \d'B^+_3)\cap C^1(\overline{B^+_3})$ satisfies
\begin{equation*}
\begin{cases}
Lu=0,\:\:u>0,&\text{in}\:B^+_3\backslash\d 'B^+_3,
\\
\a^{nb}(x)\d_bu=q(x)u,&\text{on}\:\d 'B^+_3,
\end{cases}
\end{equation*}
then there exists $C=C(n,\Lambda)>1$ such that
$$
\max_{\overline{B_1^+}} u\leq C\min_{\overline{B_1^+}}u\,.
$$
\end{lemma}

Next we will perform some computations. 
\begin{lemma}\label{int:partes}
We have:
\\\\
(a) $\int_0^{\infty}\frac{s^{\alpha}ds}{(1+s^2)^m}=\frac{2m}{\alpha +1}\int_0^{\infty}\frac{s^{\alpha+2}ds}{(1+s^2)^{m+1}}\;$, for $\alpha +1<2m$;
\\\\
(b) $\int_0^{\infty}\frac{s^{\alpha}ds}{(1+s^2)^m}=\frac{2m}{2m-\alpha -1}\int_0^{\infty}\frac{s^{\alpha}ds}{(1+s^2)^{m+1}}\;$, for $\alpha +1<2m$;
\\\\
(c) $\int_0^{\infty}\frac{s^{\alpha}ds}{(1+s^2)^m}=\frac{2m-\alpha-3}{\alpha +1}\int_0^{\infty}\frac{s^{\alpha+2}ds}{(1+s^2)^{m}}\;$, for $\alpha +3<2m$.
\end{lemma}
\begin{proof}
Integrating by parts,
$$\int_0^{\infty}\frac{s^{\alpha+2}ds}{(1+s^2)^{m+1}}=\int_0^{\infty} s^{\alpha+1}\frac{s\,ds}{(1+s^2)^{m+1}}=
\frac{\alpha+1}{2m}\int_0^{\infty}\frac{s^{\alpha} ds}{(1+s^2)^{m}}\,,$$
for $\alpha+1<2m$, which proves the item (a).

The item (b) follows from the item (a) and from
$$\int_0^{\infty}\frac{s^{\alpha}ds}{(1+s^2)^{m}}=\int_0^{\infty}\frac{s^{\alpha}(1+s^2)}{(1+s^2)^{m+1}}ds
=\int_0^{\infty}\frac{s^{\alpha}ds}{(1+s^2)^{m+1}}+\int_0^{\infty}\frac{s^{\alpha+2}ds}{(1+s^2)^{m+1}}\,.$$

To prove the item (c), observe that, by the item (a), 
$$\int_0^{\infty}\frac{s^{\alpha}ds}{(1+s^2)^{m-1}}
=\frac{2(m-1)}{\alpha+1}\int_0^{\infty}\frac{s^{\alpha+2}ds}{(1+s^2)^{m}}\,,$$ 
for $\alpha+3<2m$.
But, by the item (b), we have  
$$\int_0^{\infty}\frac{s^{\alpha}ds}{(1+s^2)^{m-1}}
=\frac{2(m-1)}{2(m-1)-\alpha-1}\int_0^{\infty}\frac{s^{\alpha}ds}{(1+s^2)^{m}}\,.$$
\end{proof}
\begin{corollary}\label{int:s}

We set $I=\int_{0}^{\infty}\frac{s^{n}}{(s^2+1)^{n}}ds\,$. Then
\\
\\
(i)$\int_{0}^{\infty}\frac{s^2+(t^2-1)}{(s^2+(t+1)^2)^{n+1}}s^{n}ds
=I\left\{\frac{n+1}{2n}(t+1)^{1-n}+\frac{n-1}{2n}(t^2-1)(t+1)^{-1-n}\right\}$;
\\
(ii)$\int_{0}^{\infty}\frac{s^2+(t^2-1)}{(s^2+(t+1)^2)^n}s^{n-2}ds
=I\left\{(t+1)^{1-n}+(t^2-1)(t+1)^{-1-n}\right\}$;
\\
(iii)$\int_{0}^{\infty}\frac{s^2+(t^2-1)}{(s^2+(t+1)^2)^{n-1}}s^{n-2}ds
=I\left\{2\frac{n-1}{n-3}(t+1)^{3-n}+2(t^2-1)(t+1)^{1-n}\right\}$.
\end{corollary}
\bp
By a change of variables we obtain
\\
\\
$\int_{0}^{\infty}\frac{s^2+(t^2-1)}{(s^2+(t+1)^2)^{n+1}}s^{n}ds
=(t+1)^{1-n}\int_{0}^{\infty}\frac{s^{n+2}}{(s^2+1)^{n+1}}ds
+(t^2-1)(t+1)^{-1-n}\int_{0}^{\infty}\frac{s^{n}}{(s^2+1)^{n+1}}ds$,
\\
$\int_{0}^{\infty}\frac{s^2+(t^2-1)}{(s^2+(t+1)^2)^n}s^{n-2}ds
=(t+1)^{1-n}\int_{0}^{\infty}\frac{s^{n}}{(s^2+1)^{n}}ds
+(t^2-1)(t+1)^{-1-n}\int_{0}^{\infty}\frac{s^{n-2}}{(s^2+1)^{n}}ds$,
\\
$\int_{0}^{\infty}\frac{s^2+(t^2-1)}{(s^2+(t+1)^2)^{n-1}}s^{n-2}ds
=(t+1)^{3-n}\int_{0}^{\infty}\frac{s^{n}}{(s^2+1)^{n-1}}ds
+(t^2-1)(t+1)^{1-n}\int_{0}^{\infty}\frac{s^{n-2}}{(s^2+1)^{n-1}}ds$.
\\\\
Then we use Lemma \ref{int:partes} to see that 
$\int_{0}^{\infty}\frac{s^{n+2}}{(s^2+1)^{n+1}}=\frac{n+1}{2n}I$, 
$\int_{0}^{\infty}\frac{s^{n}}{(s^2+1)^{n+1}}=\frac{n-1}{2n}I$,
$\int_{0}^{\infty}\frac{s^{n-2}}{(s^2+1)^{n}}=I$,
$\int_{0}^{\infty}\frac{s^{n}}{(s^2+1)^{n-1}}=2\frac{n-1}{n-3}I$ and 
$\int_{0}^{\infty}\frac{s^{n-2}}{(s^2+1)^{n-1}}=2I$.
\ep

\bigskip
\begin{lemma}\label{rel:int:t}
For $m>k+1$,
$$
\int_{0}^{\infty}\frac{t^k}{(1+t)^m}dt
=\frac{k!}{(m-1)(m-2)...(m-1-k)}\,.
$$
\end{lemma}
\begin{proof}
Integrating by parts,
$$
\int_{0}^{\infty}t^{k-1}(1+t)^{1-m}dt
=\frac{m-1}{k}\int_{0}^{\infty}t^{k}(1+t)^{-m}dt\,.
$$
On the other hand,
$$
\int_{0}^{\infty}t^{k-1}(1+t)^{1-m}dt
=\int_{0}^{\infty}\frac{t^{k-1}(1+t)}{(1+t)^m}dt
=\int_{0}^{\infty}\frac{t^{k}}{(1+t)^m}dt
+\int_{0}^{\infty}\frac{t^{k-1}}{(1+t)^m}dt\,.
$$
Hence, 
$$
\int_{0}^{\infty}\frac{t^k}{(1+t)^{m}}dt
=\frac{k}{m-1-k}\int_{0}^{\infty}\frac{t^{k-1}}{(1+t)^{m}}dt\,.
$$
Now the result follows observing that
$
\int_{0}^{\infty}\frac{1}{(1+t)^{m}}dt=\frac{1}{m-1}\,.
$
\end{proof}


\bigskip\noindent
\small{INSTITUTO DE MATEM\'{A}TICA \\UNIVERSIDADE FEDERAL FLUMINENSE \\NITER\'{O}I - RJ, BRAZIL\\E-mail addresses: {\bf{almaraz@vm.uff.br}}, {\bf{almaraz@impa.br}}

\end{document}